\documentclass[sn-aps]{sn-jnl}

\usepackage[utf8]{inputenc}
\usepackage{graphicx, xcolor}
\usepackage{tabularx}
\usepackage{amsmath,amssymb,amsfonts}
\usepackage{mathrsfs} 
\usepackage{gensymb}
\usepackage{cleveref}
\usepackage{mathtools} 
\usepackage{gensymb} 
\usepackage{caption} 
\usepackage{numprint} 
\usepackage{multirow} 
\usepackage{multicol} 
\usepackage{subfig}
\usepackage{longtable}

\usepackage{amsthm}%
\usepackage[title]{appendix}%
\usepackage{textcomp}%
\usepackage{manyfoot}%
\usepackage{booktabs}%
\usepackage{algorithm}%
\usepackage{algorithmicx}%
\usepackage{algpseudocode}%
\usepackage{listings}%

\usepackage{siunitx}
\newcommand{\sib}[1]{[\si{#1}]}

\graphicspath{{./}}

\theoremstyle{thmstyleone}%
\theoremstyle{thmstyletwo}%

\theoremstyle{thmstylethree}%

\begin{document}

\title[Article Title]{Application of Deep Learning Reduced-Order Modeling for Single-Phase Flow in Faulted Porous Media}


\author*[1]{\fnm{Enrico} \sur{Ballini}}\email{enrico.ballini@polimi.it}

\author[1]{\fnm{Luca} \sur{Formaggia}}\email{luca.formaggia@polimi.it}

\author[1]{\fnm{Alessio} \sur{Fumagalli}}\email{alessio.fumagalli@polimi.it}

\author[1]{\fnm{Anna} \sur{Scotti}}\email{anna.scotti@polimi.it}

\author[1]{\fnm{Paolo} \sur{Zunino}}\email{paolo.zunino@polimi.it}

\affil*[1]{MOX, Department of Mathematics, Politecnico di Milano, Piazza Leonardo da Vinci 32, 20133 Milano, Italy}


\abstract{We apply reduced-order modeling (ROM) techniques to single-phase flow in faulted porous media, accounting for changing rock properties and fault geometry variations using a radial basis function mesh deformation method. This approach benefits from a mixed-dimensional framework that effectively manages the resulting non-conforming mesh. To streamline complex and repetitive calculations such as sensitivity analysis and solution of inverse problems, we utilize the Deep Learning Reduced Order Model (DL-ROM). This non-intrusive neural network-based technique is evaluated against the traditional Proper Orthogonal Decomposition (POD) method across various scenarios, demonstrating DL-ROM's capacity to expedite complex analyses with promising accuracy and efficiency.}

\keywords{Porous media, faults, reduced order modeling, proper orthogonal decomposition, deep learning}


\maketitle

%

\section{Introduction} \label{sec:Introduction}
Nowadays, the practice of injecting fluids into the subsoil is gaining prominence not only for the production of fossil fuels, but also for the storage of carbon dioxide and for the strategic storage of thermal energy and the exploitation of geothermal resources. These applications are essential in our quest for sustainable and renewable energy solutions \cite{Norbotten2012}. 
Injection of fluid into the subsoil alters the local equilibrium, changing the stress field and possibly causing fault reactivation. Furthermore, leakage phenomena from storage reservoirs must be predicted \cite{Lu2012}. 
These phenomena have to be assessed taking into account the uncertainties related to the subsoil, both those concerning the physical property of the rock, e.g. the porosity or the permeability, and uncertainties about the geometry, such as the exact position of the fault or the relative displacement of the two sides. 

Due to this lack of knowledge, we need to execute multi-query applications such as sensitivity analysis, parameter identification, or uncertainty quantification, which require many evaluations of the discrete model for different scenarios. The evaluation of each scenario can be computationally demanding; therefore, these analyses can be a burden. 
Reduced order modelling (ROM) techniques come into play to provide a surrogate model that is both reliable and fast to evaluate, making such intensive computation feasible.
Reduced order models can be created using linear reduction techniques \cite{Quarteroni2016, Hesthaven2016, Brunton2019, Benner2005} with methods such as Proper Orthogonal Decomposition (POD) \cite{Quarteroni2016, Hesthaven2016, Brunton2019}, greedy algorithm \cite{Quarteroni2016, Hesthaven2016}, empirical interpolation method \cite{Quarteroni2016, Maday2009},  dynamic mode decomposition \cite{Schmid2010, Kalur2023}. Moreover, the methods can be adopted in a multi-fidelity context \cite{Peherstorfer2018} or in dynamic adaptation \cite{Peherstorfer2015}.
The primary purpose of employing data-driven ROM empowered by neural networks, as an alternative to traditional approaches, is the built-in non-intrusiveness and the ability to overcome the limitation of using a linear combination of basis functions by introducing a nonlinear trial manifold. This is achieved with the use of autoencoders, which offer advantages in highly nonlinear or advection-dominated problems \cite{Gonzalez2018, Murata2019, Hasegawa2020, Fresca2021, Fresca2021a, Fresca2022, Fresca2023, Franco2022, Franco2023}. It is also possible to discover the latent dynamic and compute a quantity of interest from it \cite{Regazzoni2023}.

In this paper, we focus on the problem of incompressible, non-reactive, single-phase flow in a faulted porous medium, considering both the physical and geometrical variability of the data. For the latter, we choose to account for the geometrical changes of the fault configuration by deforming the computational grid using a method based on algebraic equations.
We apply a new data-driven model order reduction technique based on deep feedforward neural networks, called the Deep Learning Reduced Order Model (DL-ROM) \cite{Fresca2021, Franco2022}. This method is intrinsically non-intrusive and naturally capable of efficiently dealing with nonaffine parameterizations, such as the one used for the changing geometry of the faults.
We compare DL-ROM with the well-established POD method \cite{Quarteroni2016, Hesthaven2016}, using several test cases on the problem of flows in fractured porous media with deformable geometry.

We organise the development and assessment of the proposed methodology as follows. In \Cref{sec:Mathematichal model} we present the mathematical model that governs single-phase flow in a porous medium with particular attention to the treatment of the coupling of subdomains of different dimensions. The discretization of this model is described in \Cref{sec:discretization}. \Cref{sec:rom} regards the reduced-order modelling techniques, with an introduction of both the POD and DL-ROM methods and a discussion of the nonaffine parameterisation of our problem. In \Cref{sec:mesh_deformation}, we present our methodology for deforming the geometry while maintaining a non-conforming mesh at the interfaces of subdomains. In \Cref{sec:case_studies}, we set up three different test cases to assess the main features that define the properties of the methods, such as offline and online time and solution error compared to the full order model. Finally, in \Cref{sec:multi-query} we show an example of two possible applications of the proposed model order reduction technique.

%

\section{Mathematical model} \label{sec:Mathematichal model}
We present the mathematical model of an incompressible, non-reactive, single-phase flow in a porous medium, showing the treatment of the faults. The reader can find a complete list of symbols used in this document in \Cref{appendix:nomenclature}.

\subsection{The continuous model}\label{subsec:flow_in_a_coninuous_porous_media}
Let us consider $\Omega \subset \mathbb{R}^D$, with $D = 2$ or $3$, a sufficiently regular domain that represents a porous medium, with outer boundary $\partial \Omega$ with outward normal $\upsilon$. We assume that there exists a partition of $\partial \Omega$ into two measurable parts $\partial_p \Omega$ and $\partial_q \Omega$, such that $\overline{\partial \Omega} = \overline{\partial_p \Omega} \cup \overline{\partial_q \Omega}$ and $\mathring{\partial_p \Omega} \cap \mathring{\partial_q \Omega} = \emptyset$, with $|\partial_p \Omega| \neq 0$.
We consider a Darcy model for flow in a saturated porous medium, where the Darcy velocity $q$, in \sib{\meter\per\second}, and the fluid pressure $p$, in \sib{\pascal}, satisfy the following system of partial differential equations consisting of the Darcy law and the mass balance \cite{Wangen2009}, with associated boundary conditions.
\begin{align}\label{eq:darcy}
    \begin{cases}
        q + K\nabla p =  0\\
        \nabla \cdot q = f
    \end{cases}
    & \quad \text{in } \Omega, \qquad
    \begin{cases}
        p = \overline{p} &\text{on } \partial_p \Omega, \\
        q \cdot \upsilon = \overline{q} &\text{on } \partial_q \Omega,
    \end{cases}
\end{align}
The data and parameters of the model are the permeability of the rock matrix scaled by the dynamic viscosity $K$, in $\sib{\meter^3\second\per\kilo\gram}$, a scalar forcing term $f$, in \sib{\per\second} (representing a source or a sink), and the boundary data $\overline{p}$ and $\overline{q}$,
in \sib{\pascal} and \sib{\meter\per\second}, respectively. We assume that $K$ is a bounded, symmetric, positive-definite tensor.

In this work, we consider the primal formulation associated with \eqref{eq:darcy}, where the only variable is
the pressure $p$, and the Darcy velocity may be reconstructed using Darcy's law. It reads
\begin{align}
        -\nabla\cdot(K\nabla p) =f 
    & \quad \text{in } \Omega, \qquad
     \begin{cases}
        p = \overline{p} &\text{on } \partial_p \Omega, \\
        K\frac{\partial p}{\partial \upsilon} = -\overline{q} &\text{on } \partial_q \Omega,
    \end{cases}   
\end{align}
where $\frac{\partial p}{\partial \upsilon}=\nabla p\cdot\upsilon$.
Under suitable regularity assumptions on the domain shape, forcing, and boundary terms, it is well known that the problem is well posed and admits a unique weak solution $p$ in $H^1(\Omega)$ as long as $\partial_p \Omega$ has a non-null measure, see \cite{Boffi2013}. If $|\partial_p \Omega|=0$, $p$ is known up to a constant.

\subsection{Flow in faulted porous medium}\label{subsec:flow_in_fractures}

Fractures or faults might be present in our domain of interest; these are regions with one dimension, the aperture, that may be orders of magnitude smaller than the lateral extension. For this reason, we adopt a model reduction strategy that considers fractures as objects of codimension 1.
Consequently,  we need to write a model for the flow inside the fracture by appropriately averaging the Darcy model in the equation \eqref{eq:darcy}. In the reduced model, the aperture becomes a parameter. The model applies to fractures and faults, yet because of the target application, we will use only the term faults in the sequel.

Faults might have permeabilities very different compared to the
surrounding porous media. When they are more permeable than the surrounding medium, they form a preferential flow path. 
In contrast, they block the flow and may create
compartments in the porous medium when they have a lower permeability. We are interested in both cases: for this
reason, we consider the mixed-dimensional model already used, among others, in \cite{Martin2005,
DAngelo2011,
Flemisch2016a,
Boon2018, Berre2020a}, which accounts for both situations.

\begin{figure}[h!]
    \centering
    \includegraphics[width=0.5\textwidth]{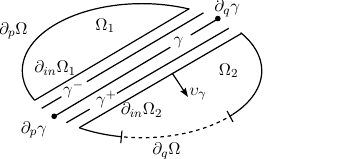}
    \caption{Generic domain divided in $\Omega_1$ and $\Omega_2$ by a fault, $\gamma$, with normal $\upsilon_\gamma$. The coupling fluxes, $\lambda^+$ and $\lambda^-$, are defined at the additional interfaces, $\gamma^+$ and $\gamma^-$. The parts of the external boundary with the Neumann conditions are indicated by $\partial_q \Omega$ and $\partial_q \gamma$, while the external boundaries with the Dirichlet conditions are denoted by $\partial_p \Omega$ and $\partial_p \gamma$. The internal boundaries facing the fault are called $\partial_{in}\Omega_1$ and $\partial_{in}\Omega_2$.}
    \label{fig:generic_domain}
\end{figure}

For the convenience of the readers, we introduce the mathematical model when only one fault, $\gamma$, is present, as shown in Fig.~\ref{fig:generic_domain}; the intersection between faults will be treated later. Moreover, we assume that the fault is planar with a unitary normal $\upsilon_\gamma$ and an
aperture $\epsilon>0$. 
The variables in the
porous medium, $\Omega$, are the same as before, while in the fault we consider the following reduced variables,
\begin{gather*}
    q_\gamma(x) = \int_{\epsilon(x)} (I - \upsilon_\gamma \otimes \upsilon_\gamma) q
    \quad \text{and} \quad
    p_\gamma(x) = \frac{1}{\epsilon(x)} \int_{\epsilon(x)} p.
\end{gather*}
Here, $q_\gamma$ is an integrated tangential flux in \sib{\square\meter\per\second}, and $p_\gamma$ an averaged pressure in \sib{\pascal}.

The normal $\upsilon_\gamma$ allows us to uniquely define the positive and negative side of the fault, as shown in Fig.~\ref{fig:generic_domain}. On each side, we introduce an additional interface, called $\gamma^+$ or $\gamma^-$, where we define a new variable $\lambda^+$ or $\lambda^-$ that represents the flux exchange between the fault and the porous medium, both measured in \sib{\meter\per\second}. We will simplify the notation by indicating $\lambda = (\lambda^+, \lambda^-)$. 

We note that now the domain $\Omega$ does not include the fault, which is, in fact, an internal boundary. Thus, the boundary of $\Omega$ can be partitioned into two open sets: $\partial_{ex}\Omega$ and $\partial_{in}\Omega$. The first indicates the external boundary, $\partial_{ex}\Omega\cap\Gamma=\emptyset$, and the second indicates the part of the boundary of $\Omega$ facing the fault. Similarly, the boundary of $\gamma$ is subdivided into two (possibly empty) disjoint parts: $\partial_{ex} \gamma$ and $\partial_{in} \gamma$, the former being the portion of the boundary of $\gamma$ in contact with $\overline{\partial_{ex}\Omega}$, 
while $\partial_{in} \gamma$ is the boundary of $\gamma$ internal to $\Omega$, characterized by $\overline{\partial \Omega} \cap \overline{\partial_{in} \gamma} = \emptyset$. The external boundary $\partial_{ex} \gamma$ is again divided into two disjoint parts, possibly empty, $\partial_{p} \gamma$ and $\partial_{q} \gamma$ where we impose pressure and fluxes, respectively. Finally, we define $\hat{\upsilon}_\gamma$ as the unit normal outward of $\partial \gamma$ (tangent to $\gamma$).

According to the mixed-dimensional model described in the cited literature, we can write the primal formulation
for the coupled problem where the unknowns are $p$ in $\Omega$, $p_\gamma$ in $\gamma$ and $(\lambda^+, \lambda^-)$ in $\gamma^+\times \gamma^-$. We have

\begin{subequations}\label{eq:full_system}
\begin{align}\label{eq:full_system_a}
    \begin{gathered}
        \begin{aligned}
            & -\nabla\cdot {K}\nabla p=f,
        \end{aligned}
        \quad \text{in } \Omega, \qquad
        \begin{aligned}
            & -\nabla \cdot \epsilon K_{\tau} \nabla p_{\gamma} + \lambda^- - \lambda^+ = f_\gamma, \quad \text{in } \gamma, \\
        \end{aligned}
    \end{gathered}
\end{align}
\begin{align}\label{eq:full_system_b}
    \begin{gathered}
        \begin{cases}
            \epsilon \lambda^+ + 2 K_n  \left( p \bigr |_{\gamma^+} - p_{\gamma} \right) = 0,
        \quad \text{on } \gamma^+, \\
        \epsilon \lambda^- + 2 K_n  \left(p_{\gamma} - p \bigr |_{\gamma^-}\right) = 0,
        \quad \text{on } \gamma^-,
        \end{cases}
    \end{gathered}
\end{align}
where \eqref{eq:full_system_a} represents the mass balances in, respectively, $\Omega$ and $\gamma$, whereas \eqref{eq:full_system_b} is the constitutive law for the fluxes $\lambda$. \color{black} $K_\tau$ and $K_n$ are the in-plane and normal permeability of the fault, scaled by the dynamic viscosity, respectively, both expressed in $\sib{\meter^3\second\per\kilo\gram}$. Furthermore, $f_\gamma$ is the averaged source or sink term in the fault, in \sib{\per\second}. Velocities can be reconstructed by Darcy's law, which for the fracture is expressed as $q_{\gamma} =- \epsilon K_{\tau} \nabla p_{\gamma}$.
Note that the differential operators in $\gamma$ are defined with respect to a local coordinate system, but we retain the same notation for simplicity. Here, with $|_{\gamma^+}$ and $|_{\gamma^-}$, we indicate the trace operators on the respective sides of the fault.

System \eqref{eq:full_system_a} with \eqref{eq:full_system_b} is complemented by the following boundary conditions
\begin{align}\label{eq:full_system_c}
    \begin{gathered}
    \begin{aligned}
            &p = \overline{p} &&\text{on } \partial_p \Omega,
            \quad
            &&K \frac{\partial p}{\partial\upsilon} = -\overline{q}, &&\text{on } \partial_q \Omega,\\
            &p_\gamma = \overline{p_\gamma} &&\text{on } \partial_p \gamma,
            \quad
            &&\epsilon K_\tau\frac{\partial p_\gamma}{\partial \upsilon_\gamma} = -\overline{q_\gamma} &&\text{on } \partial_q \gamma,
        \end{aligned}\\
        \frac{\partial p_\gamma}{\partial \upsilon_\gamma}=0 \quad \text{on } \partial_{in} \gamma.
    \end{gathered}
\end{align}
\end{subequations}
The last relation, known as the tip condition, imposes a null flux.

\subsection{Intersections}
\begin{figure}[h!]
    \centering
    \includegraphics[width=0.4\textwidth]{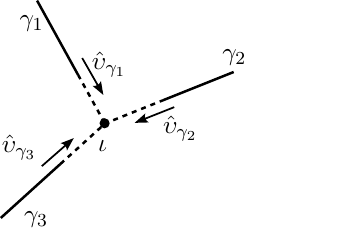}
    \caption{Y-shaped intersection. The three faults, $\gamma_1$, $\gamma_2$, $\gamma_3$, cross each other at the point $\iota$. $\hat{\upsilon}_{\gamma_i}$ are the outward unit vectors aligned with the respective fault.}
    \label{fig:intersecton}
\end{figure}

We focus here on the intersection of 1D faults because it is the only kind of intersection that we consider in this paper. Different intersection topologies can exist, such as X-, T-, or Y-shaped intersections; see Fig.~\ref{fig:intersecton} for an example of a Y-shaped intersection. All types of intersection are treated similarly, as explained below. The intersection point is represented as a 0D domain where the mass conservation equation and the coupling conditions must be solved. The 0D mass conservation is: 
\begin{equation*}
    \sum_{i=0}^{n_{br}} \lambda_{\gamma_i} = 0,
\end{equation*}
where $n_{br}$ is the number of intersecting branches, $\lambda_{\gamma_i}$ $\sib{\meter\per\second}$ are the fluxes exchanged between the branches through the intersection point, $\iota$: $\lambda_{\gamma_i} = q_{\gamma_i}|_\iota \cdot \hat{\upsilon}_{\gamma_i}$, where $\hat{\upsilon}_{\gamma_i}$ is the outwards unit vector aligned with the $i$-th fault.
Regarding the coupling condition, we can use a model considering pressure jumps as:
\begin{equation*}
    \epsilon^2 \lambda_\gamma + 2 K_\iota \left(p_\iota-p_\gamma|_\iota \right) = 0,
\end{equation*}
where $K_\iota$ $\sib{\meter^3 \second \per \kilo\gram}$ is a representative value of permeability at the
intersection, $p_\iota$ $\sib\pascal$ is the pressure of the intersection.
The well-posedness of the problem of flows in faulted porous media, including intersections, for a different, yet equivalent, formulation has been studied in \cite{Formaggia1998, Nordbotten2018}.

\section{Discretization}\label{sec:discretization}

We discretize problem \eqref{eq:full_system} with a cell-centered finite volume method \cite{Blazek2015, Hirish2007}. The degrees of freedom of the primary variables $p$, $p_\gamma$, and $\lambda$ are located in the center of the cells that discretize the domain and fractures.
Numerical fluxes at the cell boundary are computed with the multipoint flux approximation (MPFA) \cite{Aavatsmark2002, Nordbotten2021, Starnoni2019, Stefansson2018}, in particular, the so-called MPFA-O method first proposed in~\cite{Aavatsmark2002} implemented in the Porepy library \cite{Keilegavlen2019}. 
We stress that our objective is to solve the problem for different values of the geometrical and physical parameters, respectively $\mu_{geom}$ and $\mu_{phy}$. We set $\mu = (\mu_{geom}, \mu_{phy}) \in \Theta \subset \mathbb{R}^e$, where $\Theta$ is the space of parameters, and $e$ the total number of parameters. 
Variations in geometric parameters lead to domain deformations.  
We detail the issues related to domain discretization in \Cref{sec:mesh_deformation}.


For a given grid, the discretization of the full-order model \eqref{eq:full_system_a} leads to a perturbed saddle point linear system of the form
\begin{equation}\label{eq:linear_sys_N}
     \begin{bmatrix}
        A & B_1^T \\
        B_2 & -C
    \end{bmatrix} 
    \begin{bmatrix}
        \hat{\mathfrak{p}} \\
        \boldsymbol{\lambda}
    \end{bmatrix}= b_N,
\end{equation}
where $A = \begin{bmatrix}
    A_p & 0 \\
    0 & A_{p_\gamma}
\end{bmatrix}$ is a semi-definite positive matrix containing the discretization of the fluxes in \eqref{eq:full_system_a}, precisely $A_p$ discretizes the fluxes in $\Omega$ and $A_{p_\gamma}$ discretizes the fluxes in $\gamma$. The terms
$B_1$ and $B_2$ derive from the coupling conditions, involving the fluxes $\lambda$, while $C$ is the mass matrix that appears in the discretization of the constitutive law \eqref{eq:full_system_b}.
We note that \eqref{eq:linear_sys_N} satisfies the assumptions for well-posedness; for example, see \cite{Boffi2013} Sect. 5.5.
The unknown $\hat{\mathfrak{p}} = [\mathfrak{p}^T, \mathfrak{p}_\gamma^T]^T$ contains the degrees of freedom $\mathfrak{p}$ and $\mathfrak{p}_\gamma$ related to $p$ and $p_\gamma$, while $\boldsymbol{\lambda}$ contains the degrees of freedom of $\lambda$. The previous system can be rewritten in the following compact form.
\begin{equation}\label{eq:linear_sys_N}
    A_N u_N = b_N.
\end{equation}
Here, $A_N \in \mathbb{R}^{N\times N}$, $u_N \in \mathcal{V}_N = \text{col}(A_N)\subset \mathbb{R}^{N}$ is the unknown vector consisting of all degrees of freedom $u_N = [\hat{\mathfrak{p}}^T, \boldsymbol{\lambda}^T]^T$, $b_N\in \mathbb{R}^{N}$, $N$ is the number of degrees of freedom of the full order model. 
All these quantities depend on $\mu$, which is not explicitly indicated to simplify the notation.

%

\section{Model order reduction}\label{sec:rom}

For several applications, it is necessary to repeatedly query a model with various input values to perform tasks such as sensitivity analysis, parameter estimation through inverse problem solving, or uncertainty quantification. When the model is complex and computationally demanding, managing multiple queries becomes impractical. Therefore, the development of a fast and reliable surrogate model is crucial.

We address this problem by approximating the \textit{solution manifold}, $\mathcal{S} := \{u_N\}_{\mu \in \Theta}$, of the flow problem in faulted porous media \eqref{eq:full_system} through a reduced model that takes as input the physical and geometric parameters $\mu_{phy}$, $\mu_{geom}$, and returns an estimate of the corresponding variables $p$, $p_\gamma$, $\lambda$. The key idea under the model order reduction procedure is that a few main patterns characterize the parametric dependence of a problem described by many degrees of freedom (d.o.f.). Thus, we exploit these patterns to formulate a smaller problem that is easier to solve. Then,  we reconstruct an approximation of the desired solution using a linear or
nonlinear map $\mathcal{M}$ from the reduced space to the full-order space.

Model order reduction can be model-based or data-driven \cite{Quarteroni2016, Brunton2019, Benner2015}. The former uses the model equations in their differential or discrete form to derive a reduced version. It has the property of being intrusive, which means that we need to manipulate the equations and, consequently, modify the solver code to implement the reduced model. 
In the latter case, instead, the reduced model is built using a collection of data through which the dynamics of the system can be inferred; the procedure is non-intrusive because it is equations and software-agnostic.

In the following sections, we apply two reduced-order model techniques to the problem of flows in a porous medium with faults. The first is the Proper Orthogonal Decomposition (POD), which is model-based and linear. The second technique, called DL-ROM, is data-driven and approximates the nonlinear map $\mathcal{M}$ using neural networks.


\subsection{Proper orthogonal decomposition}\label{sec:pod}
In this case, the map $\mathcal{M}$ is linear, and it is represented by an orthogonal matrix $\Phi \in \mathbb{R}^{N \times n}$ called \textit{transition matrix}, with $n < N$. A slightly different definition is considered in \Cref{sec:dl_rom}, where a non-linear map replaces the linear map, $\Phi$. The matrix $\Phi$ contains information about the basis functions that represent the full-order solutions. Typical ways to compute it are the singular value decomposition (SVD) and the greedy algorithm \cite{Hesthaven2016, Quarteroni2016}. We now briefly recall the SVD process.
A number $n_s$, $n < n_s < N$, of full-order model solutions, called snapshots, are calculated and collected in the snapshot matrix,
\begin{equation}\label{eq:snapshot_matrix}
    S =
    \begin{bmatrix}
        u^{(1)} | \ldots | u^{(n_s)}
    \end{bmatrix}.
\end{equation}
Then, the singular value decomposition factorizes $S$ as $S = U \Sigma V^H$,
%
%
where $U \in \mathbb{R}^{N \times N}$ is a unitary matrix whose columns are the left singular vectors, $\Sigma \in \mathbb{R}^{N \times n_s}$, is a rectangular diagonal matrix containing the singular values in decreasing order, $V \in \mathbb{R}^{n_s \times n_s}$ is a unitary matrix whose columns are the right singular vectors.
$\Phi$ will be formed by the first $n$ columns of $U$ associated with the highest singular values $\Phi = U_{tr}$,
where $[U_{tr}]_{ij} = [U]_{ij}$, $i = 1, \ldots, N$, $j = 1, \ldots, n$. This choice of $\Phi$ is optimal in the sense that it minimizes the reconstruction error in the snapshots:
\begin{equation*}
    \sum_{i=1}^{n_s} \| u^{(i)}-\Phi \Phi^\top u^{(i)}\|_2^2 = \min_{W \in \mathbb{R}^{N\times n}} \sum_{i=1}^{n_s} \| u^{(i)} - W W^\top u^{(i)} \|_2^2.
\end{equation*}
The proof follows the Eckart-Young theorem \cite{Eckart1936}, originally discovered by Schmidt in the continuous framework \cite{Schmidt1907}; see \cite{Quarteroni2016} for a detailed proof. The previous equation tells us that the POD is an optimal linear reduced-order modeling approach in the sense described previously. 

Given the basis functions stored in the columns of the transition matrix $\Phi$,
we define the \textit{reconstructed solution}, $\Tilde{u}_N$, belonging to the full-order space $\mathcal{V}_N$, obtained from the reduced basis solution:
\begin{equation}\label{eq:tilde{u}_N_u_n}
    \Tilde{u}_N = \Phi u_n.
\end{equation}
where $u_n \in \mathcal{V}_n \subset \mathbb{R}^n$ are the \textit{reduced coefficients}, solution of the reduced problem. 

A possible strategy to retrieve $u_n$ is to make the residual orthogonal to the linear subspace defined by the column of
$\Phi$, i.e., to compute a Galerkin projection into the
linear subspace $S_n = \mathrm{col}(\Phi)$,
\begin{equation*}
    \Phi^\top \left( b_N - A_N \Phi u_n \right) = 0,
\end{equation*}
which yields
\begin{equation*}
    \Phi^\top A_N \Phi u_n = \Phi^\top b_N.
\end{equation*}
Defining $A_n = \Phi^\top A_N \Phi u_n$, and $b_n = \Phi^\top b_N$, we have the following reduced problem:
\begin{equation}\label{eq:reduced_system}
    A_n u_n = b_n.
\end{equation}
After the reduced problem has been solved, $\Tilde{u}_N$ is retrieved from \eqref{eq:tilde{u}_N_u_n}.

Since our problem features multiple coupled variables, an alternative formulation involves applying the SVD separately to each variable. In this way, we replace the transition matrix with a block-diagonal one. Considering a problem with two generic physical variables $c$ and $d$ we have:
\begin{align*}
    S_c &=
    \begin{bmatrix}
        u_c^{(1)} | \ldots | u_c^{(n_s)}
    \end{bmatrix} = U^c \Sigma^c (V^c)^H \\
    S_d &=
    \begin{bmatrix}
        u_d^{(1)} | \ldots | u_d^{(n_s)}
    \end{bmatrix} = U^d \Sigma^d (V^d)^H.
\end{align*}
%
In addition, the truncation is done independently, producing the matrices $U^c_{tr}$ and $U^d_{tr}$ whose elements are: $U^c_{tr,ij} = U^c_{ij}$, $i=1,\ldots,N$, $j=1,\ldots,n_c$ and $U^d_{tr,ik} = U^d_{ik}$, $i=1,\ldots,N$, $k=1,\ldots,n_d$.
The transition matrix is a block diagonal matrix of the truncated left singular vectors matrices:
\begin{equation*}
    \Phi =
    \begin{bmatrix}
    U_{tr}^{c} & 0 \\
    0 & U_{tr}^{d}
    \end{bmatrix}.
\end{equation*}
We refer to this approach as block-POD. Considering the unknowns of our problem, $p$, $p_\gamma$, $\lambda$, the matrix takes the form:
\begin{equation*}
    \Phi = 
    \begin{bmatrix}
    U_{tr}^{p} & 0 & 0 \\
    0 & U_{tr}^{p_\gamma} & 0 \\
    0 & 0 & U_{tr}^{\lambda}
    \end{bmatrix}.
\end{equation*}
Here, we choose the number of modes in each left-singular vector to be the same for each physical variable $n_{p} = n_{p_\gamma} = n_{\lambda}$.

A problem is affine if the discrete operator can be rewritten as a linear combination where only the coefficient depends on $\mu$. Taking into account the full-order matrix, $A_N$, the affine parametric dependence implies that:
\begin{equation*}\label{eq:affine}
    A_N = \sum_{q=1}^Q \theta_q(\mu) A_N^q,
\end{equation*}
where $\theta_q(\mu)$ are scalar functions, and $A_N^q$ are constant matrices, so they are computed once and for all, independently of $\mu$. The reduced order matrix, $A_n$, is thus expressed by the following sum:
\begin{equation}\label{eq:affine_reduced}
    A_n = \Phi^\top \sum_{q=1}^Q \theta_q(\mu) A_N^q \Phi = \sum_{q=1}^Q \theta_q(\mu) \Phi^\top A_N^q \Phi = \sum_{q=1}^Q \theta_q(\mu) A_n^q.
\end{equation}
POD is particularly efficient when the affine parametric dependence is satisfied because the calculation of $A_n$ for each new value of the parameters $\mu$ does not require reassembling the matrix $A_N$  but using \eqref{eq:affine_reduced} with the precomputed $A_n^q$.
However, as we will see in \Cref{sec:mesh_deformation} the changes in geometry related to the sliding of a fault imply a non-affine problem.



\subsection{Deep Learning-ROM}\label{sec:dl_rom}
We restrict ourselves to the core ideas of the DL-ROM approach, and the reader is invited to read \cite{Fresca2021, Fresca2022, Fresca2021a, Fresca2023, Franco2022, Franco2023} for more details. 
In this case, the parameter-to-solution map $\mathcal{M}$ is approximated by employing a neural network that naturally accounts for the nonlinear structure of the solution manifold.
The reduction procedure is divided into two parts.

\begin{figure}
    \centering
    \includegraphics[width=0.7\textwidth]{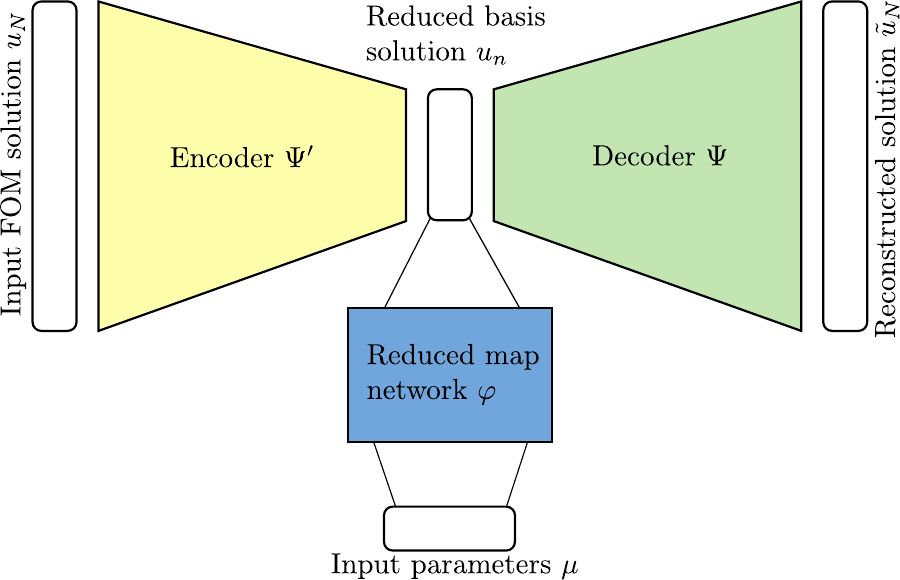}
    \caption{DL-ROM scheme. The encoder, $\Psi'$, takes the FOM solution, $u_N$, and returns the reduced basis solution, $u_n$. The decoder $\Psi$ reproduces an approximate solution, called reconstructed solution $\Tilde{u}_N$, with the only knowledge of the reduced solution. The reduced map network, $\varphi$, approximates $u_n$ as a function of the given parameter $\mu$. }
    \label{fig:dl-rom scheme}
\end{figure}

\begin{description}
\item (i) First, we perform a dimensionality reduction step that identifies a low-dimensional latent space obtained from a nonlinear mapping of the full-order space to $\mathbb{R}^n$. With little abuse of notation, here $n$ denotes the dimension of the reduced space, which can be different from the corresponding space obtained using the POD method. The mappings from the full-order space to the reduced-order space and back, named $\Psi': \mathcal{S} \rightarrow \mathcal{V}_n$ and $\Psi: \mathcal{V}_n
\rightarrow \mathcal{V}_N$ respectively,
are approximated here by feedforward neural networks. Here, $\mathcal{V}_n$ denotes the reduced-order space obtained with the encoder. More precisely, the approximation of $\Psi'$ is called the \emph{ encoder} neural network, while $\Psi$ is called the \emph{ decoder}. Their combination is called the \emph{ autoencoder}; see Fig.~\ref{fig:dl-rom scheme} for a graphical representation of the adopted neural networks.
This approach defines a nonlinear trial manifold because the d.o.f. in the latent space is mapped to the full-order space by the nonlinear function provided by the neural network. Note that neural networks work in discrete spaces; in particular, the encoder input is the discrete solution $u_N$ on a given mesh, and the decoder outputs the reconstructed solution on the same mesh.

\item (ii) Once the reduced space has been identified, we have to surrogate the operator that solves the problem in the reduced space. This is done again using a data-driven approach. Specifically, a third neural network, named $\varphi$, is trained to provide a representation of the solution in the reduced space for any admissible value of the parameters. This approach is feasible due to the low dimensionality of the input and output spaces. In essence, we employ a deep feedforward neural network called \textit{reduced map nerwork}, to describe the map, $\varphi$, from the parameter space to the reduced solution space: $\varphi: \Theta \rightarrow \mathcal{V}_n$. Therefore, the reduced solution is $u_n = \varphi(\mu)$. 
\end{description}

In conclusion, combining steps (i) and (ii) of this procedure, the resulting approximation of the parameter-to-solution map is the following composition: $\Psi(\varphi(\cdot))$, so the reconstructed solution is
\begin{equation*}\label{eq:sol_recon_dl}
    \Tilde{u}_N = \Psi(u_n) = \Psi(\varphi(\mu)).
\end{equation*}

In \cite{Franco2022}, it has been observed that this method enjoys some optimality properties in terms of the reducibility of the solution manifold. Given the manifold $\mathcal{S}$, the authors of \cite{DeVore1989} define the nonlinear counterpart of Kolmogorov $n$-width, $\delta_n(\mathcal{S})$, as the worst reconstruction error using the best encoder and decoder maps $\Psi$ and $\Psi'$, respectively:
\begin{equation*}
    \delta_n(\mathcal{S}) = \inf_{\scriptsize{\begin{matrix}
        \Psi' \in \mathcal{C}(\mathcal{S}, \mathcal{V}_n) \\
        \Psi \in \mathcal{C}(\mathcal{V}_n, \mathcal{V}_N)
    \end{matrix}}} \sup_{u \in \mathcal{S}} \| u - \Psi (\Psi'(u)) \|.
\end{equation*}
Then, following \cite{Franco2022}, the \textit{minimal latent dimension} of $\mathcal{S}$, denoted as $n_{\mathrm{min}}(\mathcal{S})$, is defined as the smallest $n$ for which $\delta_n(\mathcal{S}) = 0$.
Under the hypothesis that the map $\mu \rightarrow u_N$ is continuous and injective and $\Theta$ has a non-empty interior, it is shown in \cite{Franco2022} that $n_{\min}(\mathcal{S}) = e$. In other words, there exists an autoencoder, with the bottleneck width equal to the size of the parameter space, capable of achieving a null reconstruction error.


The encoder and decoder process the discrete solution as a singular vector containing all physical components $p$, $p_\gamma$, $\lambda$.

\subsubsection{Training}
Training the neural network can follow the stages reported in \Cref{sec:dl_rom}.
In the dimensionality reduction step, the encoder and decoder networks, $\Psi'$ and $\Psi$, must be trained.
We associate with this task the following loss function:
\begin{equation}\label{eq:loss_function_1}
    \ell_1 = \frac{1}{N} \| u_N(\mu_i) - \Psi(\Psi'(u_N(\mu_i))) \|_2^2.
\end{equation}
This is an unsupervised learning task in which the autoencoder $\Psi \circ \Psi'$ is trained using a sequence of unlabeled data representing randomly sampled points in the solution manifold.
The second step is to train the reduced map network to learn the map $\varphi$. This is done by minimizing the following loss function:
\begin{equation}\label{eq:loss_function_2}
    \ell_2 = \frac{1}{n} \| \Psi'(u_N(\mu_i)) - \varphi(\mu_i) \|_2^2.
\end{equation}
In this case, we use a sequence of labeled data $[\mu_i,\Psi'(u_N(\mu_i)]$, representing the parameter-to-reduced state map at a suitable number of points in the parameter space.

Two options are available to train the entire network $\Phi \circ \varphi$. One is to proceed as previously described, minimizing the function $\ell_1$ first and then the function $\ell_2$. This method has the advantage of reducing the computational complexity of training because the parameters of the networks $\Psi'$ and $\Psi$ are optimized independently of those of $\varphi$.
The second training approach consists of optimizing all networks together, by minimizing the total loss function defined as
\begin{equation}\label{eq:loss_function}
    \mathscr{L} = \alpha \ell_1 + \beta \ell_2,
\end{equation}
where $\alpha >0 \in \mathbb{R}$, $\beta >0 \in \mathbb{R}$, are user-defined hyper-parameters.
We remark that the encoder is used only during the offline stage, see \Cref{sec:offline_online}. The query of the network for a new value of the parameter, $\mu$, does not require the evaluation of the encoder. Indeed, the final reduced model is composed of only the reduced map network and the decoder.
It appears that two-stage training is not beneficial, as suggested by \cite{Fresca2021}. Training the three neural networks together is advantageous in terms of both accuracy and speed.

\subsection{Offline-Online stages}\label{sec:offline_online}

Model order reduction strategies are typically divided into two stages: the \textit{offline} phase and the \textit{online} phase. The offline phase encompasses all operations that do not need to be repeated when querying the reduced model with a new parameter value, denoted as $\mu$. This phase includes data generation, which is common to both the POD and DL-ROM methods. Specifically, for POD, this phase also involves calculating $\Phi$ and, in cases of affine parametric dependence, computing $A_n^q$. For DL-ROM, the offline phase ends with neural network training.

Although the offline phase can be time intensive, it is a one-time process. Once completed, the online phase begins. For POD, this includes the assembly of $A_n$, solving the reduced system, and reconstructing the solution. On the contrary, the online phase of DL-ROM simply involves the forward evaluation of neural networks $\varphi$ and $\Psi$.

%

\section{Mesh deformation} \label{sec:mesh_deformation}
We change the domain by deforming a reference computational grid to account for geometric uncertainties. This strategy has the advantage of keeping the same grid topology and having a fixed number of degrees of freedom. As a consequence,  all snapshots have the same size. Moreover, the ordering of the d.o.f. is kept, so they remain associated with the same physical region since we consider small geometric deformations. The general procedure presented in the following can also work for large deformations in both the 2D and 3D cases at the price of additional mesh handling complexity. 
We adopt a technique based on radial basis functions (RBF) to deform the mesh, and we show first the standard approach \cite{Boer2007, Forti2014, Aubert2017}
, then its adaptation to our specific problem of geometries containing sliding faults.  

In the standard approach, we seek to interpolate each component of the displacement function, $s(x) \in \mathbb{R}^{D}$, $D$ being the physical dimension, using a linear combination of a given radial basis function, $g(d)$, where $d = d(x_1, x_2) = \|x_1-x_2\|$ is the distance between the two points. By defining $g^*(x_1, x_2) = g(d(x_1,x_2))$, we have the following.
\begin{equation}\label{eq:rbf_basic}
    s(x) = \sum_{j=1}^l g^*(x, x_{c_j}) \zeta_j,
\end{equation}
where $\zeta_j \in \mathbb{R}^{D}$, are the unknown coefficients and  $x_{c_j}$ are the so-called control points. We use the following radial basis function: $g(d) = d/0.2$.
To find $\zeta_j$, we need to apply the constraints formed by selecting a number $l$ of relevant points of which we have information about their displacement. For instance, we can choose the nodes on the fixed borders of the domain or the nodes belonging to elements whose geometry is uncertain. These points are called control points, $x_{c_j}, \ j=1, \ldots, l$, and that is where we enforce a known displacement, $\overline{s}_j$:
\begin{equation}\label{eq:constr_displ}
    s(x_{c_j}) = \overline{s}_j.
\end{equation}
Evaluating \eqref{eq:rbf_basic} at the control points, we have a linear system for each component, $m = 1,\ldots,D$, of the vector displacement, $s$:
\begin{equation}\label{eq:rbf_sys_basic}
    G z_m = \sigma_m,
\end{equation}
where $G \in \mathbb{R}^{l \times l}, \ [G]_{i,j} = g(d(x_{c_i}, x_{c_j}))$, $z_m \in \mathbb{R}^l$ is the unknown vector where
$[z_m]_j = [\zeta_m]_j$. The right-hand side $\sigma_m \in \mathbb{R}^l$ contains the known displacement as a function of the geometrical parameters $[\sigma_m]_j = [\overline{s}_j]_m(\mu_{geom})$.

\begin{figure}[h]
    \centering
    \includegraphics[width=0.4\textwidth]{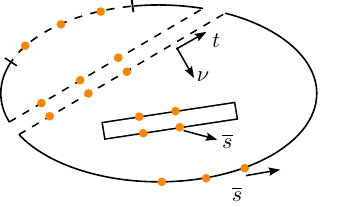}
    \caption{Generic faulted domain with control points, in orange. The dashed lines represent sliding surfaces, one is on $\partial_{ex}\Omega$, and the control points therein belong to $C_{s}$, the other two represent the boundaries $\partial_{in}\Omega$ facing a fault, the control points on those lines are included in the set $C_{sf}$. The gap between the two faces of the fault was added for graphical reasons only, and it is actually absent, so the control points on the fault may be coincident. Other control points in $C_d$ are placed on the fracture where a rigid displacement, $\overline{s}$, is applied. The remaining control points in $C_{df}$ are on the boundary of $\Omega$ where a displacement $\overline{s}$ is enforced.}
    \label{fig:control_points}
\end{figure}

We aim to let the two sides of the fault slide independently, so the control points must be positioned on both sides to avoid any undesirable deformation of the fault. Some control points inevitably become very close or even coincident, causing $G$ to be ill-conditioned or singular. Moreover, to make the process more practical, instead of the displacement constraint, we would like to impose a sliding constraint on $x_{c_j}$ on some specific surfaces, for example, the fault surface or the boundary of the domain; see Fig.~\ref{fig:control_points}. 

We now show an improvement of the standard mesh deformation method to meet our requirements. Taking a generic surface, $S_i$, which could be, for example, a fault or a boundary face, with normal $\nu_i$ and two (one in 2D) non-parallel tangent unit vectors $t_i$ and $b_i$, the sliding constraint can be described by two conditions: the non-penetration condition:
\begin{equation}\label{eq:no_penetration}
    s(x_{c_j}) \cdot \nu_i = 0,
\end{equation}
and by the no-tangential contribution condition:
\begin{align}\label{eq:no_tangential}
    \begin{gathered}
    \zeta_j \cdot t_i = 0, \\
    \zeta_j \cdot b_i = 0,
    \end{gathered}
\end{align}
where $\zeta_j$ is the coefficient associated with $x_{c_j}$ on the sliding surface. Equation \eqref{eq:no_penetration} sets the control points on a sliding surface free to move except in the normal direction. Equation \eqref{eq:no_tangential} implies that the tangential displacement of the control points on the fault does not affect the displacement of all the other points. Anyhow, the displacement of the control points on the fault is influenced by all the other control points.

To address the issue of coincident points and add the sliding condition, let us define an index set $C$ of $x_{c_j}$ and divide it into four subsets: $C_d$ containing the control points not on the fault where displacement is enforced, $C_{df}$ containing the control point on the fault where displacement is enforced, $C_s$ and $C_{sf}$ containing the sliding control points, respectively, not in the fault and in the fault. Moreover, we call $\overline{C}_s$ the set of surfaces where sliding conditions are applied. We introduce the \textit{side function} of a point, $\beta(x, \nu)$, with respect to the sliding fault, $\gamma$, with normal $\nu$:
\begin{equation*}\label{eq:side_fcn}
    \beta(x, \nu) = \frac{ \mathrm{sign}\left( (x-x_{ref})\cdot \nu \right) +1 }{2},
\end{equation*}
$x_{ref} \in \gamma $ is a reference point.
We define the influence function, $\mathcal{I}$, as a function that hides the points that are not on the same side of the fault. An expression could be (for the sake of simplicity, we restrict the discussion to the case of a single sliding fault, with normal $\nu_2$):
\begin{multline*}
\mathcal{I}(x, \beta(x, \nu_2), \beta(x_{c_j}, \nu_2), \nu_2) =\\
    \begin{cases}
        1  &\mathrm{if} \ x_{c_j} \in C_{ds}
        \\
        NXOR(\beta(x, \nu_2), \beta(x_{c_j}, \nu_2)) &\mathrm{if} \ x_{c_j} \in C_{sf} \ \mathrm{and} \ x \in \partial_{in}\Omega  \\
        \dfrac{|(x-x_{ref})\cdot \nu_2|}{\lVert (x-x_{ref}) \rVert} NXOR(\beta(x, \nu_2), \beta(x_{c_j}, \nu_2)) &\mathrm{if} \ x_{c_j} \in C_{sf} \ \mathrm{and} \ x \in \Omega,
    \end{cases}
\end{multline*}
where  $C_{ds} = C_d \cup C_s$ and $C_f = C_{df} \cup C_{sf}$  and $NXOR$ is the combination of the logical operations ``not" and ``xor", $NXOR(a,b)=\neg(a \oplus b)$.
Introducing $\mathcal{I}$ to \eqref{eq:rbf_basic} and the displacement and sliding constraints, setting $g^\dag(x, x_{c_j}, \nu_2) = \mathcal{I}(x, \beta(x, \nu_2), \beta(x_{c_j}, \nu_2), \nu_2) g^*(x, x_{c_j}, \nu_2)$ the system to be solved to deform the mesh becomes:

\begin{equation*}
\begin{cases}
    s(x) = \sum_{j=1}^l g^\dag(x, x_{c_j}, \nu_2) \zeta_j, \\
    s(x_{c_j}) \cdot \nu_i = 0, &  x_{c_j} \in C_s \cup C_{sf}, \ S_i \in \overline{C}_s \\
    \begin{aligned}
        & \zeta_j \cdot t_i = 0,   \\
        & \zeta_j \cdot b_i = 0,
    \end{aligned}
    & x_{c_j} \in C_s \cup C_{sf} \ s.t. \ x_{c_j} \  \text{on} \ S_i \in \overline{C_s} 
\end{cases}
\end{equation*}
For example, for a 2D case with one sliding surface, labelled as 1, and one sliding fault, labelled as 2, the system becomes:
\begin{equation*}
    \begin{bmatrix}
        G \quad 0 \\
        0 \quad G \\
        G \nu_{1x} \quad G \nu_{1y} \\
        G \nu_{2x} \quad G \nu_{2y} \\
        H_1 t_{1x} \quad  H_1 t_{1y} \\
        H_2 t_{2x} \quad H_2 t_{2y} \\
    \end{bmatrix}
    \begin{bmatrix}
        z_x \\
        z_y
    \end{bmatrix}
    =
    \begin{bmatrix}
        \overline{s}_x \\
        \overline{s}_y \\
        0 \\
        0 \\
        0 \\
        0
    \end{bmatrix}
\end{equation*}
with $[G]_{ij} = g^\dag(x_{c_i}, x_{c_j}, \nu_2)$ for $x_{c_i} \in C_d \cup C_{df}$ and $x_{c_j} \in C$,
and
\begin{align*}
    [H_{1}]_{ij} =
    \begin{cases}
        1 \quad \mathrm{if}\ x_{c_i}, x_{c_j} \in C_s, \\
        0 \quad \mathrm{else}.
    \end{cases}
    \quad
    [H_{2}]_{ij} =
    \begin{cases}
        1 \quad \mathrm{if}\ x_{c_i}, x_{c_j} \in C_{sf}, \\
        0 \quad \mathrm{else}.
    \end{cases}
\end{align*}
We aim now to demonstrate that the relative sliding of the two sides of a fault introduces a non-affine problem, thereby reducing the efficiency of the POD reduced order modeling technique, as discussed in \Cref{sec:pod}. From \eqref{eq:full_system}, the pressure on the internal boundaries of $\Omega$ facing the fault $\gamma$, see Fig.~\ref{fig:generic_domain}, can be rewritten as:
\begin{equation}\label{eq:p_p_equality}
    p|_{\partial_{in}\Omega_2} = p_{\gamma} + \frac{\epsilon}{2 K_n}\lambda^+.
\end{equation}
Let $r$ be the combination of the maps $r_1: \partial_{in}\Omega_2 
 \rightarrow \gamma^+$ and $r_2: \gamma^+
 \rightarrow \gamma$  map from $\partial_{in}\Omega_2$ to $\gamma$, so $r = r_1 \circ r_2$. A generic point $x_\gamma$ on the fault corresponds to $x_\gamma = r(x_{\partial_{in}\Omega_2}; \mu)$ where $r$ depends on the geometric parameterization governed by $\mu_{geom}$. The map $r$ could be non-affine, which already spoils the problem's affinity.
 Introducing $r$ in \eqref{eq:p_p_equality} to explicitly write the link between the pressures in the two subdomains, we have:
\begin{equation*}
    p(x_{\partial_{in}\Omega_2}) =
    p_{\gamma}(r(x_{\partial_{in}\Omega_2}; \mu)) + \frac{\epsilon}{2 K_n} \lambda_1(r(x_{\partial_{in}\Omega_2}; \mu)),
\end{equation*}
where the dependence on $\mu$ cannot be separated, even if $r$ were an affine map. In the weak formulation, this implies that the calculation of the pressure on the internal boundaries contains terms that cannot be separated from the parameters, so the problem is not affine.

\begin{figure}[h]
    \centering
    \subfloat[original mesh]{\includegraphics[width=4cm]{ 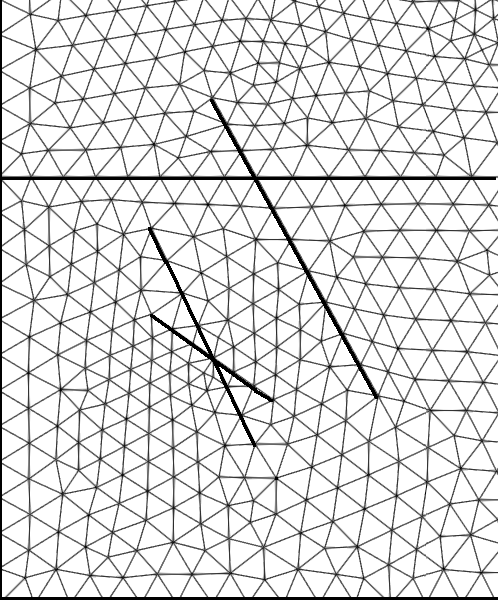}}
    \hfill
    \subfloat[deformed mesh]{\includegraphics[width=4cm]{ 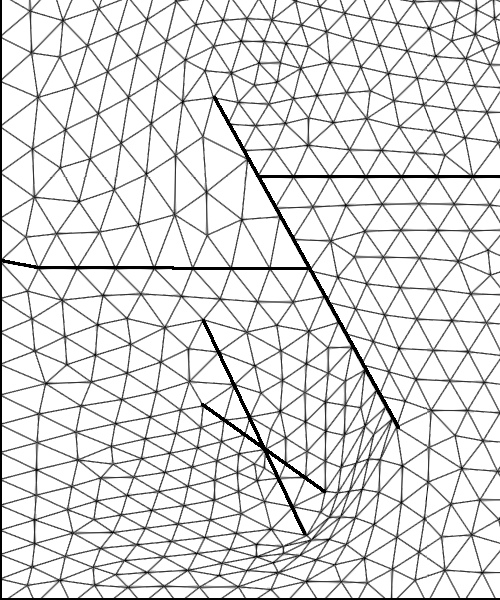}}
    \caption{Mesh deformation. Zoom of the bottom-left corner of the domain of the third test case, see \Cref{sec:test_3}. A rigid displacement is imposed to the left side of the fault and to the two small intersecting fractures. The left boundary is a sliding surface where sliding conditions are applied, on the bottom one, instead, a null displacement is enforced.}
    \label{fig:deformed_mesh}
\end{figure}

%

\section{Numerical validation}\label{sec:case_studies}
The reduced order model techniques described above have been applied to three 2D and 3D test cases with different specifics. The first test is a simple problem with a 2D domain in which uncertainties are associated with the permeability of the rock and the throw of the fault. The second case is the three-dimensional extension of the first one. The third case is a complex fracture network with varying boundary conditions.

In all the case studies we assess the quality of the reduced model by measuring the maximum, minimum, and averaged relative norm of the error, defined as:
\begin{align}\label{eq:errors}
    \begin{gathered}
        e_{max} = \max_j( \|e(\mu_j)\|_2 / \|u_N(\mu_j)\|_2 ), \\
        e_{min} = \min_j( \|e(\mu_j)\|_2 / \|u_N(\mu_j)\|_2 ), \\
        e_{ave} = \frac{1}{J}\sum_{j=1}^J \|e(\mu_j)\|_2 / \|u_N(\mu_j)\|_2,
    \end{gathered}
\end{align}
where $J$ is the size of the test dataset and $e(\mu)$ is the difference between the FOM and ROM solution, $e(\mu) = u_N(\mu)-\Tilde{u}_N(\mu)$.

\subsection{Case 1 - setup}\label{sec:test_1}   

\begin{figure}[h]
    \centering
    \includegraphics[width=0.5\textwidth]{ 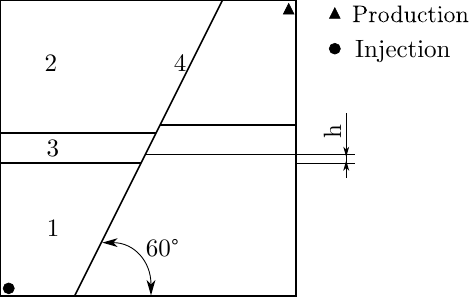}
    \caption{Case study 1. Geometry. The fault cuts the domain from the bottom to the top with an angle of $60\degree$. The right layers can slide along the fault. The production is in the top right corner and the injection is on the bottom-left corner.}
    \label{fig:test_1_domain}
\end{figure}

As a first test case, we studied a unit 2D square domain cut by a fault whose aperture is $\varepsilon = 10^{-3}$ tilted by an angle of $60\degree$ from the horizontal, Fig.~\ref{fig:test_1_domain}. There are two layers separated by a caprock, and they are identified by different values of permeabilities $K_1$, $K_2$, $K_3$.
An injection point and a production point are placed, respectively, at the bottom-left and top-right corners, where we enforce counterbalanced fluid fluxes. We impose homogeneous Neumann boundary conditions at all boundaries except for the bottom-left corner, where the pressure is set to 1.

The parameters of the problem are the three values of permeability, $K_1 \in [10^{-2}, 10^{-1}]$, $K_2  \in [10^{2}, 10^{3}]$, $K_3 \in [10^{-4}, 10^{-3}]$, of the three layers, the permeability $K_4 \in [10^{-4}, 10^{-3}]$, of the fault such that $K_\tau = K_4$, $K_n = 2K_4/\varepsilon$, and the height $h \in [0, 0.07]$, of the right horizons: a small displacement is allowed along the fault direction. The maximum displacement value $h = 0.07$ is less than the thickness of the caprock equal to 0.1.
This setting implies a pressure distribution with three main regions: high pressure layer at the bottom, low pressure layer at the top, and high-pressure gradient in the caprock that horizontally separates the domain, as we can see in Fig.~\ref{fig:test_1_single_snap}(a) that represents a typical snapshot of the test data set. Furthermore, based on the specific values of the permeabilities, there may be a jump in pressure across the two faces of the fault. 

Subsoil permeabilities may vary widely, even by several orders of magnitude, so it is convenient to express them as exponential functions and sample the exponent. Therefore, the permeabilities are written in the following form: $K_i = e^{\eta_i}, i=1,\ldots,4$ and the parameter vector consists of the exponents of the permeabilities and the height of the horizons: $\mu = (\eta_1, \eta_2, \eta_3, \eta_4, h)$.

The data set is made up of 1000 snapshots whose parameters are randomly sampled from a uniform distribution (the same strategy is adopted for the other tests). It is divided into a training dataset (800 snapshots), a validation dataset (100 snapshots), and a test dataset (100 snapshots). These data sets are used in both the POD approach (Sect.~\ref{sec:test_1_pod}) and the DL-ROM approach (see Section~\ref{sec:test_1_dlrom}).

Note that the training data set may be small compared to the complexity of the problem and the number of parameters. A comparison with a uniform grid sampling reveals that we are sampling the parameters with fewer than 4 points for each of the 5 axes of the space $\Theta$. This is because we aim to apply the reduced-order modeling techniques to real scenarios, where a single snapshot generation is extremely expensive, making the generation of very large datasets infeasible, and therefore the test cases are designed accordingly.

\subsubsection{Case 1 - POD}\label{sec:test_1_pod}
\begin{figure}[h]
    \centering
    \subfloat[POD, DL-ROM]{\includegraphics[width=5cm]{ 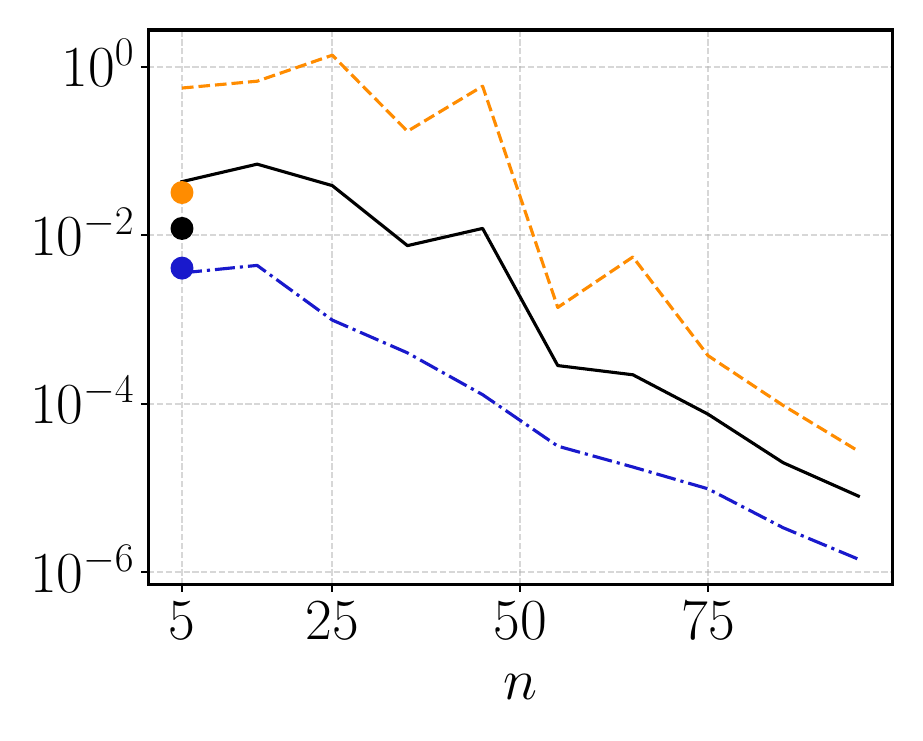}}
    \subfloat[block-POD]{\includegraphics[width=5cm]{ 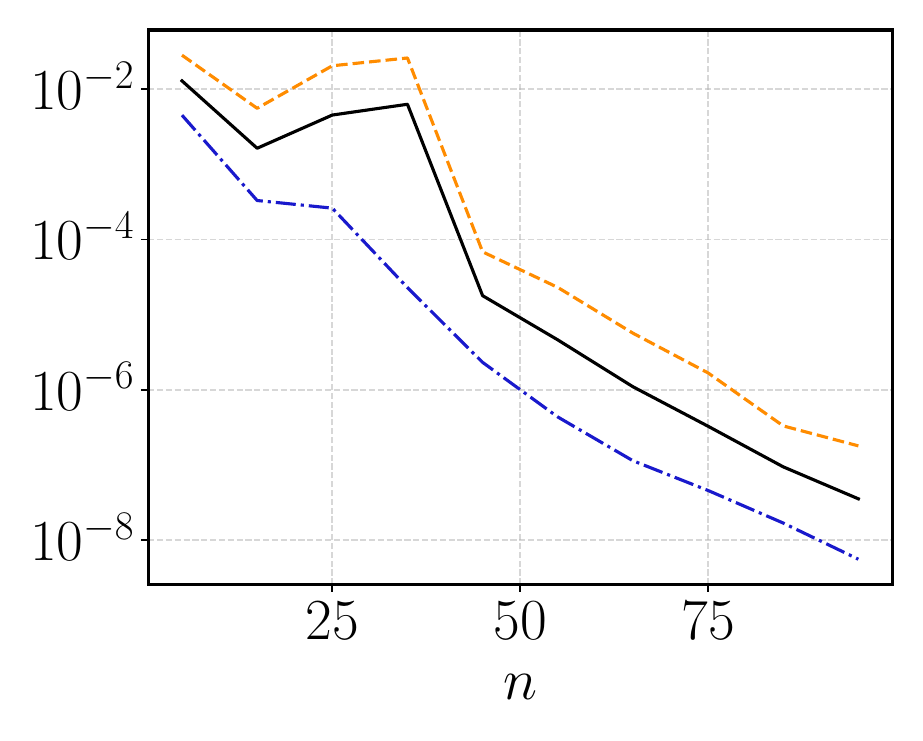}} \\
    \vspace{0.2cm}
    \includegraphics[width=10cm]{ 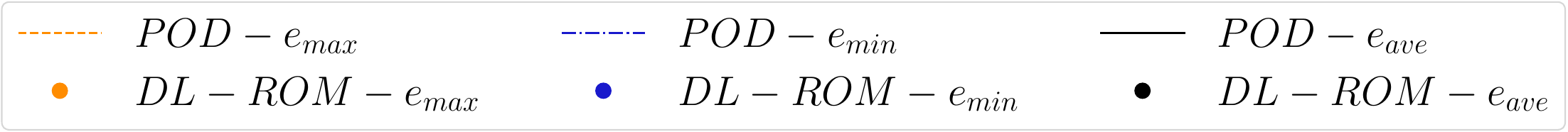}
     \caption{Case study 1. $e_{max}$, $e_{min}$, $e_{ave}$, respectively, maximum, minimum and average errors versus the size of the reduced space, $n$, of the reduced model obtained with (a) standard POD and (b) block-POD. The points on the panel (a) represents the value errors for the DL-ROM, described in \Cref{sec:test_1_dlrom}.}
     \label{fig:test_1_err}
\end{figure}

\begin{figure}[h]
    \centering
    \includegraphics[width=7cm]{ 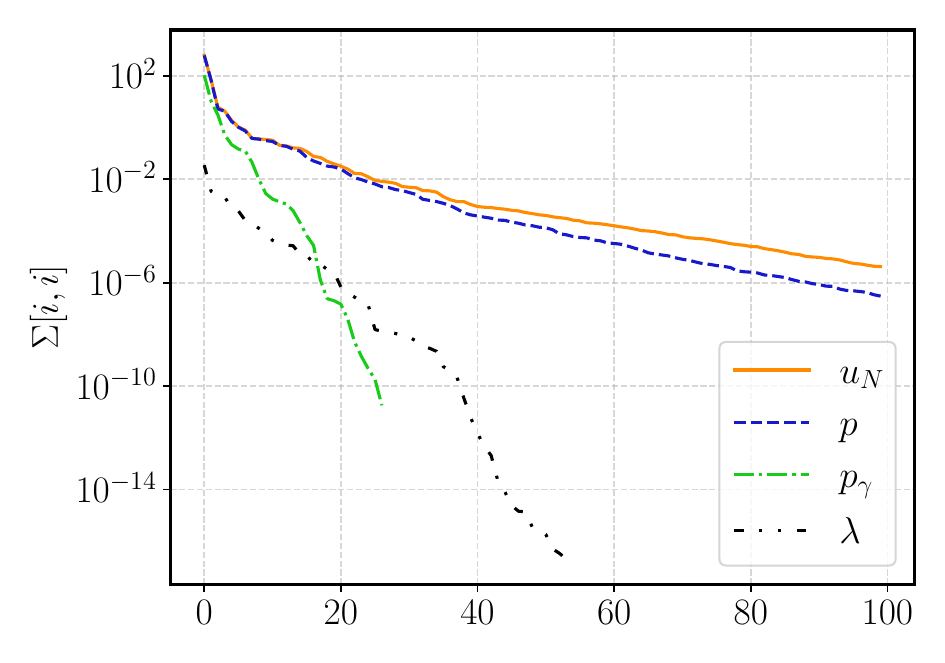}
     \caption{Case study 1. The decay of singular values of the snapshot matrix of the fully coupled POD approach for $u_N$, compared to the block-POD case applied to each variable $p_\gamma$, $\lambda$ and $p$ independently.}
     \label{fig:test_1_singular_vals_allofthem}
\end{figure}

Both POD and block-POD (\Cref{sec:pod}) are shown. After the offline phase, the quality of the reduced model is assessed using the test data set. Each parameter of the test data set is used to generate the reduced solution and then the reconstructed solution, which is compared with the solution of the full-order model to compute the errors defined in \eqref{eq:errors}. We evaluated the errors for different sizes of the reduced space, $\mathcal{V}_n$, which corresponds to different truncations of the transition matrix, as shown in Fig.~\ref{fig:test_1_err}. The block POD error shown in Fig.~\ref{fig:test_1_err}(b) decreases faster than the standard POD; we justify this trend by examining the decay of the singular values of each variable in Fig.~\ref{fig:test_1_singular_vals_allofthem}. We observe a faster decay in the singular values of the POD block snapshot matrices computed for $p_\gamma$ and $\lambda$ compared to those related to $p$. The trend of the latter resembles that of the complete discrete vector $u_N$. Then, decoupling the variables takes advantage of the faster decay of the approximations for $p_\gamma$ and $\lambda$, which otherwise would not have an effect in the fully coupled approach.

Some oscillations are observed in the trend of the error; nevertheless, it decreases to negligible levels taking a sufficient number of modes. 

Offline time encompasses the generation of the training data set and all computation related to the SVD decomposition required to obtain the transition matrix. Indicatively, data generation takes 4 min 30 s on an Intel i5-1135G7 and a fraction of a second to finalize the offline phase.
The online phase includes the construction of the matrix that assembles the linear system, its solution, and the reconstruction of the solution. Its time is assessed at 30 points of the test data set. The main statistics are shown in Fig.~\ref{fig:test_1_online_time}, as well as the times for the FOM and DL-ROM. As in the other tests, similar results are obtained with the block-POD. Indeed, in both the POD and block-POD methods, the online time is dominated by the assembly of the full order model matrix $A$, a common step for both methods. Consequently, for the sake of simplicity, the timing results of block-POD will not be reported. The square is delimited by the first and third quartiles, whereas the whiskers extend from the box by 1.5x the interquartile range. The yellow line represents the median. The number of modes $n=45$ is chosen because the related error is similar to the error obtained with the DL-ROM, so we can compare the timing for the same error. The consequent compression rate is $\frac{n}{N_h} = 4.3\%$.
%
%
%

%

%
\subsubsection{Case 1 - DL-ROM}\label{sec:test_1_dlrom}
When not specified differently, we use the following architectures and settings in each test case: the encoder is made of a fully connected network with one hidden layer, the decoder architecture is symmetric with respect to that of the encoder, and the reduced map network is made of a fully connected network with one hidden layer. As a nonlinear activation function, $\rho$, we use, for all networks, the PReLU function: $\mathrm{PReLU}(x) = \max(0,x) + a \min(0,x)$, where $a$ is a trainable weight. The weights of the layers $m$ are initialized from the uniform distribution $\mathcal{U}\left(-\sqrt{1/\mathrm{size}(m-1)}, \sqrt{1/\mathrm{size}(m-1)}\right)$.

The data used to train, validate, and test neural networks are the same as those used for the POD approach.

The input of the encoder is normalized to $[0, 1]$. Since the ratio of the physical variables may differ by some order of magnitude, we decided to normalize each discrete physical variable taking individually their maximum and minimum values over the training dataset.

Details of the neural networks of this test case are summarized in Tab~\ref{tab:case_1_dense}. With the use of neural networks, in this case, we reach a compression rate of $\frac{n}{N} = 0.5\%$, which is almost ten times higher than that of the POD. 

The network is trained for 4000 epochs with Adam optimizer~\cite{Kingma2014}, using the hyperparameters suggested in the cited reference and a minibatch size of 32. The learning rate is initialized to $10^{-3}$ and is reduced by a factor of 0.6 every 500 epochs. Once the networks are trained, the reconstruction errors are: $e_{max} = 3.2\%$, $e_{min} = 0.41\%$, and $e_{ave} = 1.2\%$, they are depicted in Fig.~\ref{fig:test_1_err} (a) as points. The distance between the maximum and minimum values is small, which means that the network is accurate evenly throughout the parameter space. 

\begin{figure}
    \centering
    \includegraphics[width=7cm]{ 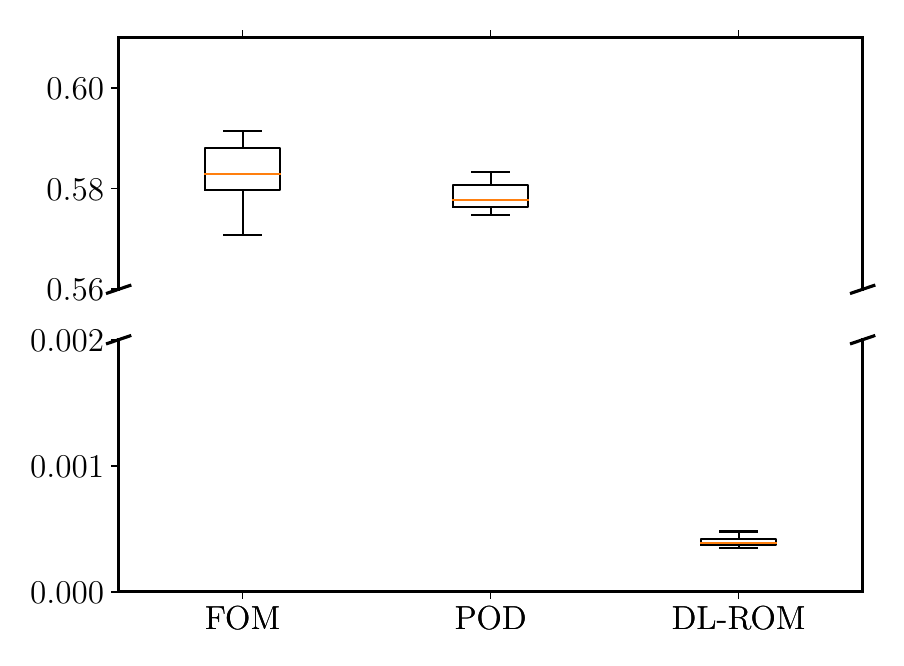}
    \caption{Case study 1. Evaluation times in seconds of the full-order model, the reduced model obtained through POD, and the reduced model obtained through the DL-ROM. Evaluations of reduced models are timed on 30 points of the test dataset.}
    \label{fig:test_1_online_time}
\end{figure}

The data generation and the training of the neural networks constitute the offline phase; the latter takes about 40 minutes to complete the 4000 epochs. The online phase is represented by the evaluation of the reduced map network and the decoder; this time it is evaluated on the Nvidia MX330 graphic card and is shown in Fig.~\ref{fig:test_1_online_time}. We can observe that, despite a slightly longer offline phase, the online is almost three orders of magnitude faster than the POD case. 
In fact, the online time savings (normalized difference of times) with respect to the FOM model is $99.14 \%$. We want to highlight that the DL-ROM allows for the efficient exploitation of both CPU and GPU hardware, while the POD does not run efficiently on GPUs.

\begin{figure}
    \centering
    \subfloat[FOM]{\includegraphics[width=4cm]{ 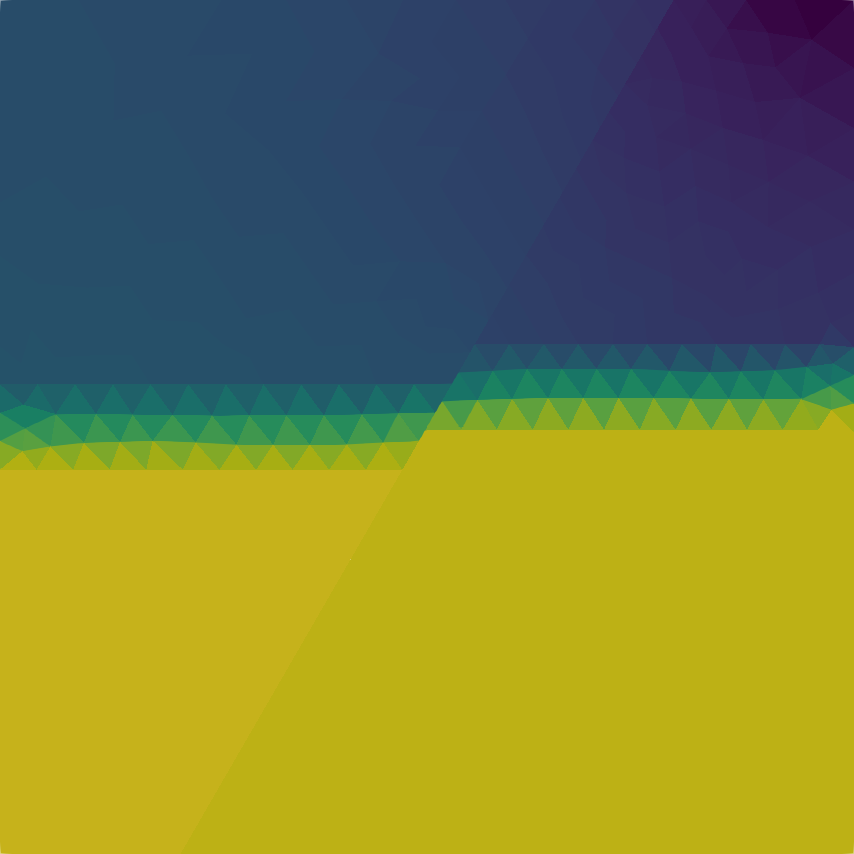}}
    \hspace*{0.1cm}
    \subfloat[POD]{\includegraphics[width=4cm]{ 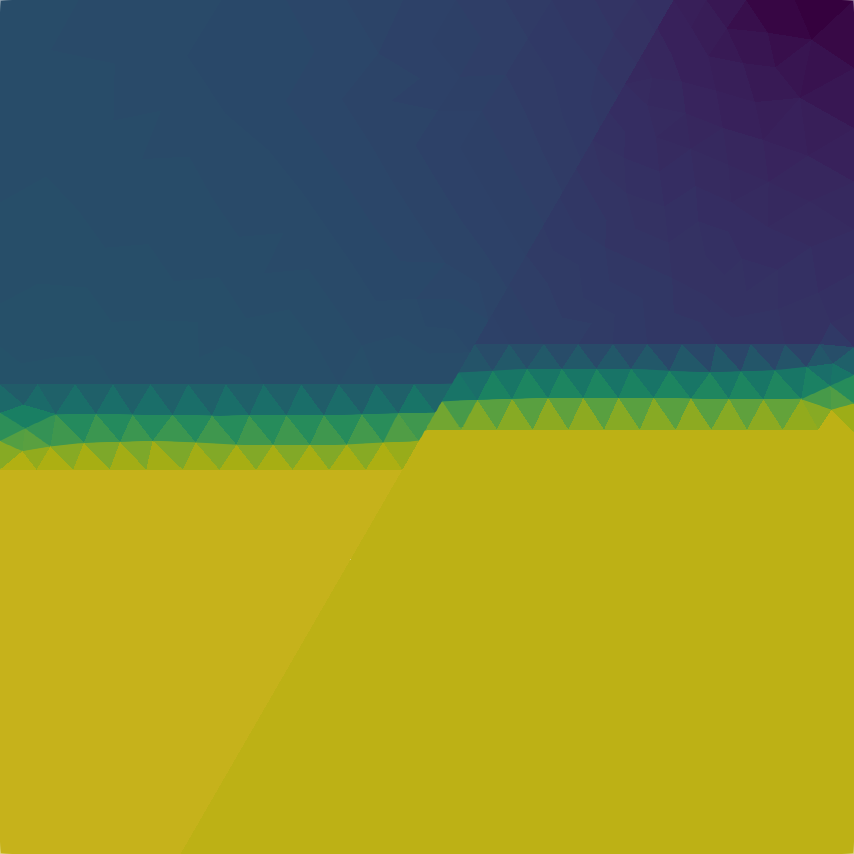}}
    \hspace*{0.1cm} \\
    \subfloat[block-POD]{\includegraphics[width=4cm]{ 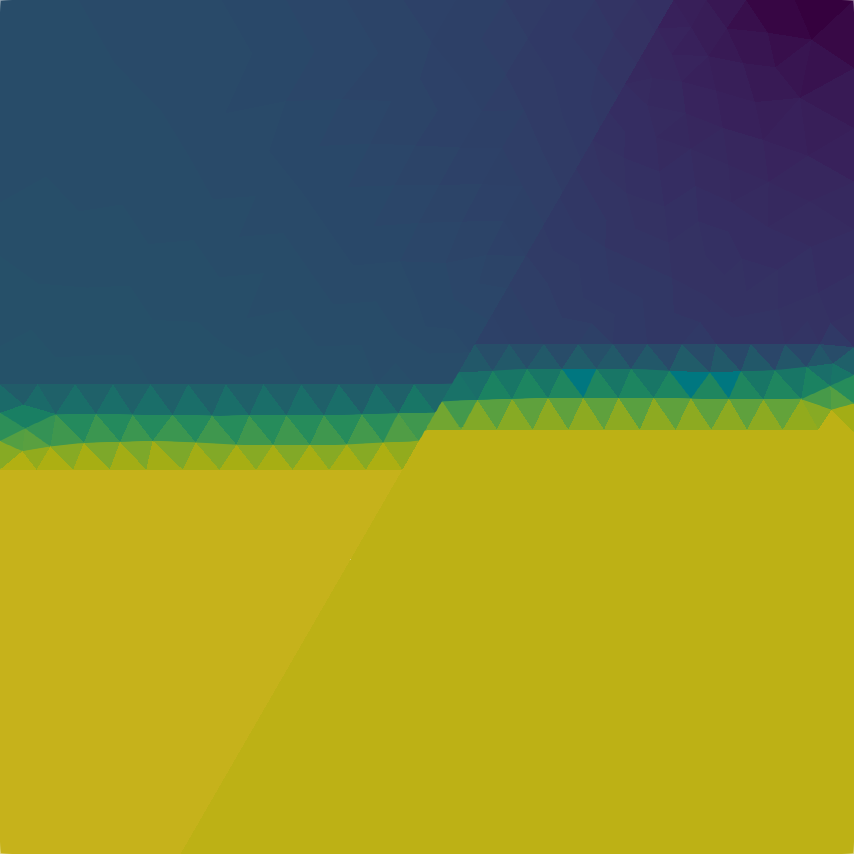}}
    \hspace*{0.1cm}
    \subfloat[DL-ROM]{\includegraphics[width=4cm]{ 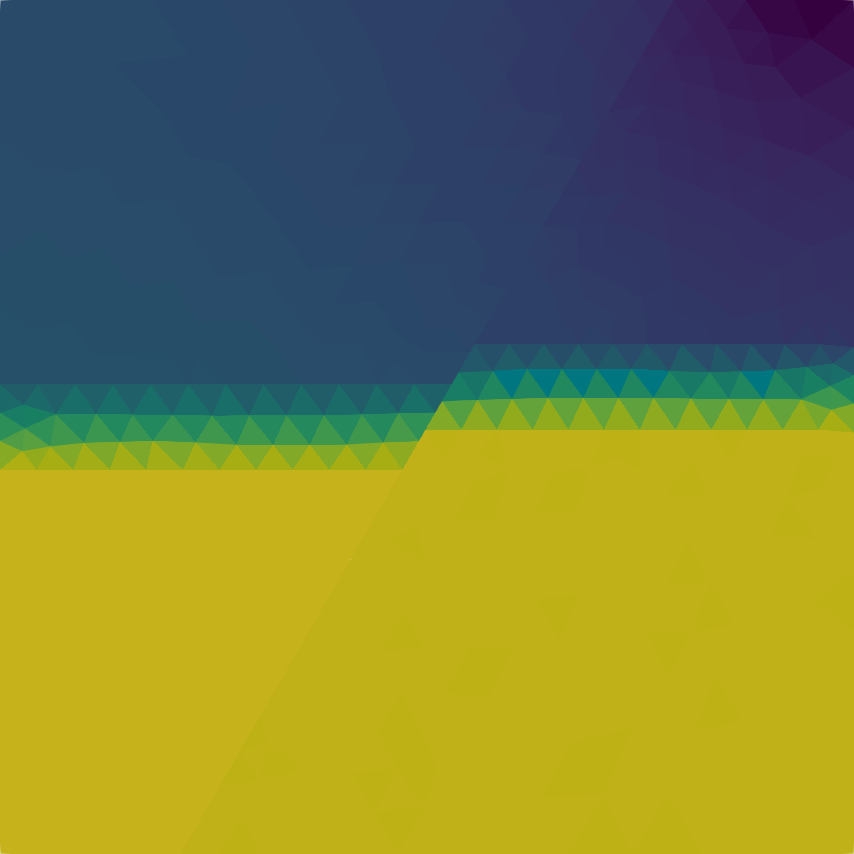}} \\
    \vspace{0.3cm}
    \includegraphics[width=4cm]{ 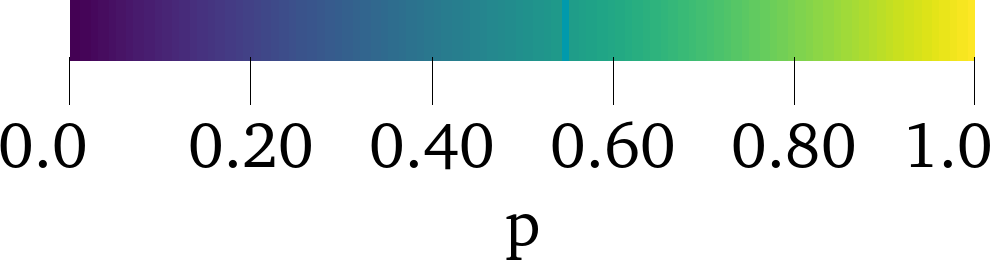}
    \caption{Case study 1. Reconstructed solutions for a specific value of the parameter $\mu$. (a) FOM. (b) POD. (c) block-POD (d) DL-ROM.
    }
    \label{fig:test_1_single_snap}
\end{figure}

A visualization of the reconstructed solution compared to the full-order model is illustrated in Fig.~\ref{fig:test_1_single_snap}. We can observe qualitatively that the reconstructed solutions of both POD and DL-ROM are close to the FOM solution in the entire domain.
\begin{table}[]
\begin{tabular}{llll}
                                     & layer        & \# tw    & \# tw tot \\
                                     \hline \\
\multirow{3}{*}{Encoder}             & fc, 1045, 181, PReLU  & \multirow{2}{*}{\numprint{190238}} & \multirow{8}{*}{\numprint{382622}}  \\
                                     & fc, 181, 5, PReLU     &                         &                          \\
                                     &                       &                         &                           \\
\multirow{3}{*}{Decoder}             & fc, 5, 181, PReLU     & \multirow{2}{*}{\numprint{191277}} &                          \\
                                     & fc, 181, 1045, PReLU  &                         &                          \\
                                     &                       &                         &                           \\
\multirow{2}{*}{Reduced map network} & fc, 5, 100, PReLU     & \multirow{2}{*}{\numprint{1107}}   &                          \\
                                     & fc, 100, 5, PReLU     &                         & \\
                                     \hline
\end{tabular}
\caption{Case study 1. Architecture of the networks that make up the reduced model. ``fc" stands for fully connected layer, tw is the number of trainable weights.}
\label{tab:case_1_dense}
\end{table}




%

\subsection{Case 2 - setup}
\begin{figure}[h]
    \centering
    \includegraphics[width=0.5\textwidth]{ 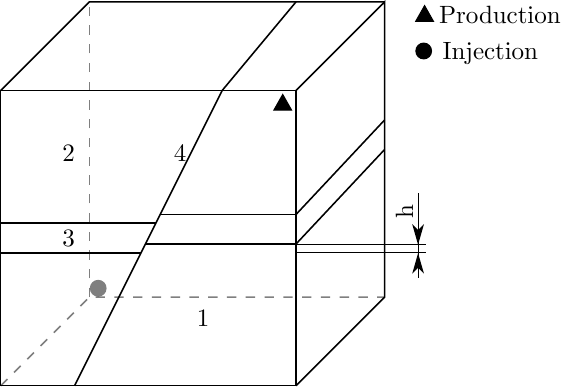}
    \caption{Case study 2. The domain of this test corresponds to the geometry of case 1 extruded along the z-direction. Production and injection are placed at the two opposite corners.}
    \label{fig:test_2_domain}
\end{figure}

The second test is a 3D simulation that resembles the first in most of its features. The main differences are the domain (Fig.~\ref{fig:test_2_domain}) which is a unit cube whose normal section to $z$ corresponds to the domain pictured in Fig.~\ref{fig:test_1_domain}, and the position of the injection and production points that are at two opposite corners of the cube. We chose this test case to study the effects of adding the third physical dimension without adding many other factors.

The datasets are made similarly to the previous test case: the parameters are sampled from a random distribution, 800 snapshots constitute the training dataset, 100 snapshots are used for the validation dataset and 100 for the test dataset.

\subsubsection{Case 2 - POD}
\begin{figure}
    \centering
    \subfloat[POD, DL-ROM]{\includegraphics[width=5cm]{ 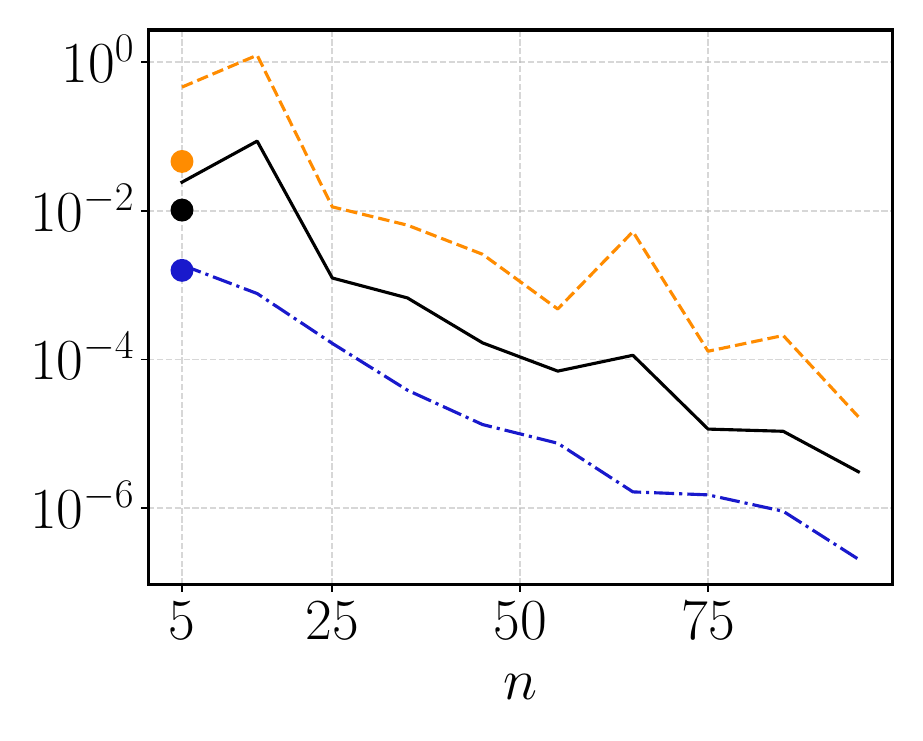}}
    \subfloat[block-POD]{\includegraphics[width=5cm]{ 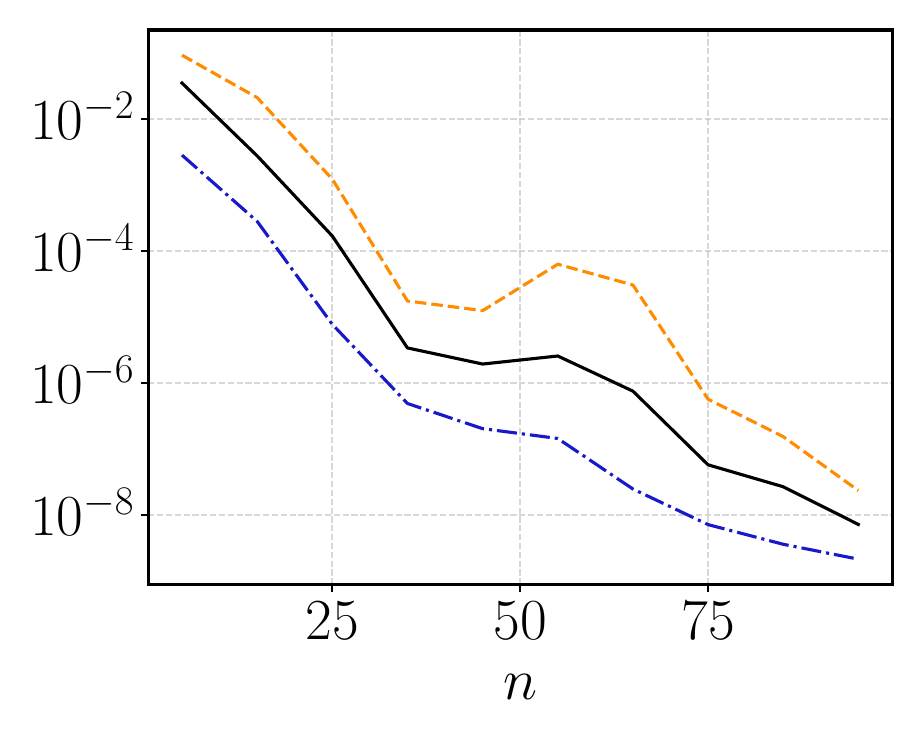}}
    \vspace{0.2cm}
    \includegraphics[width=10cm]{ test_1_legend.pdf}
    \caption{Case study 2. $e_{max}$, $e_{min}$, $e_{ave}$, respectively, maximum, minimum, and average error versus the size of the reduced space, $n$, of the reduced model obtained with (a) standard POD, (b) block-POD. The points on panel (a) represent the value errors for the DL-ROM, described in \Cref{sec:test_2_dlrom}.}
    \label{fig:test_2_err}
\end{figure}

Fig.~\ref{fig:test_2_err} shows the errors defined by \eqref{eq:errors} with respect to the size of the reduced space. We notice a slightly slower decrease of the error with respect to test case 1 and that, for this scenario, the block POD proves to have higher performance than the standard one. The online time depends on the number of modes used; to facilitate the comparison with the other methodology, we report in Fig.~\ref{fig:test_2_online_time} the online time for the standard POD for $n = 18$ which entails a test error similar to that of the DL-ROM. The related compression rate is equal to $\frac{n}{N} = 0.2\%$. Due to the great number of d.o.f., the online time reaches a high value of about 15 seconds, Fig.~\ref{fig:test_2_online_time}, it is mainly due to the computations for deforming the mesh and the rediscretization of the problem due to its non-affine nature. 
%
%
\subsubsection{Case 2 - DL-ROM}\label{sec:test_2_dlrom}
\begin{figure}
    \centering
    \subfloat[FOM]{ \includegraphics[width=0.35\textwidth]{ 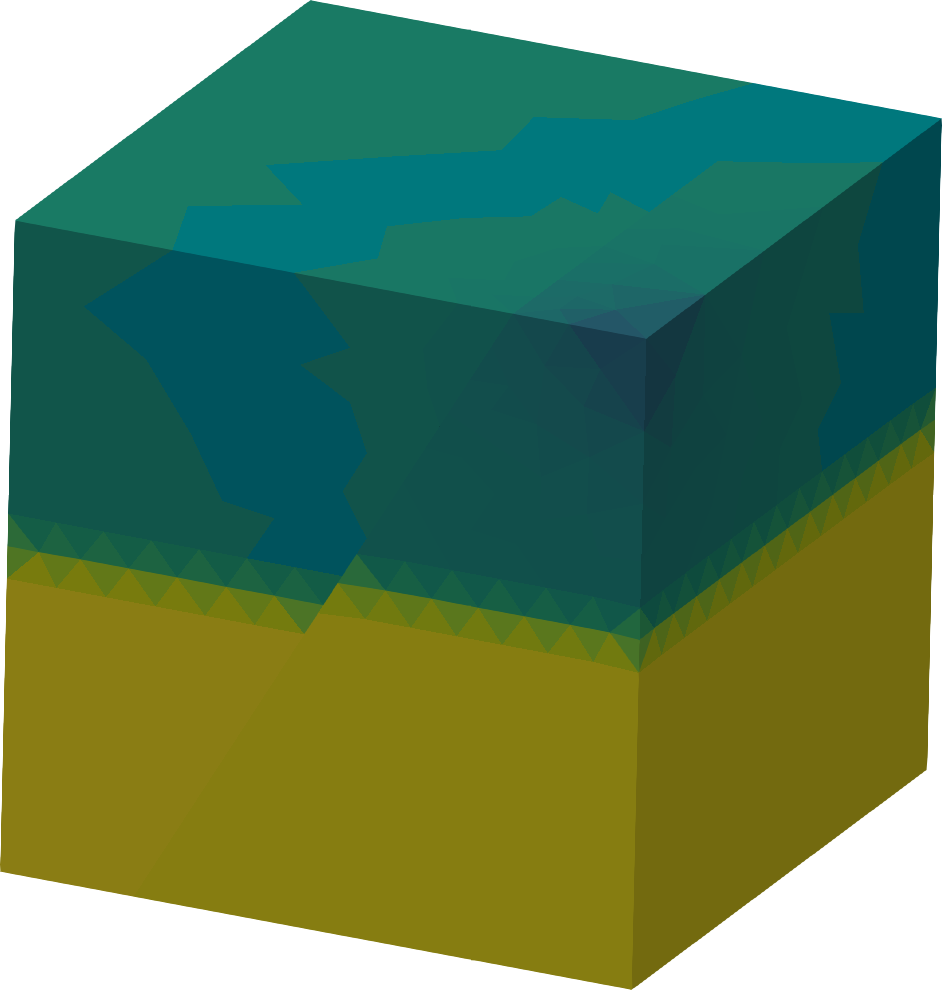} }
    \hspace*{0.1cm}
    \subfloat[POD]{\includegraphics[width=0.35\textwidth]{ 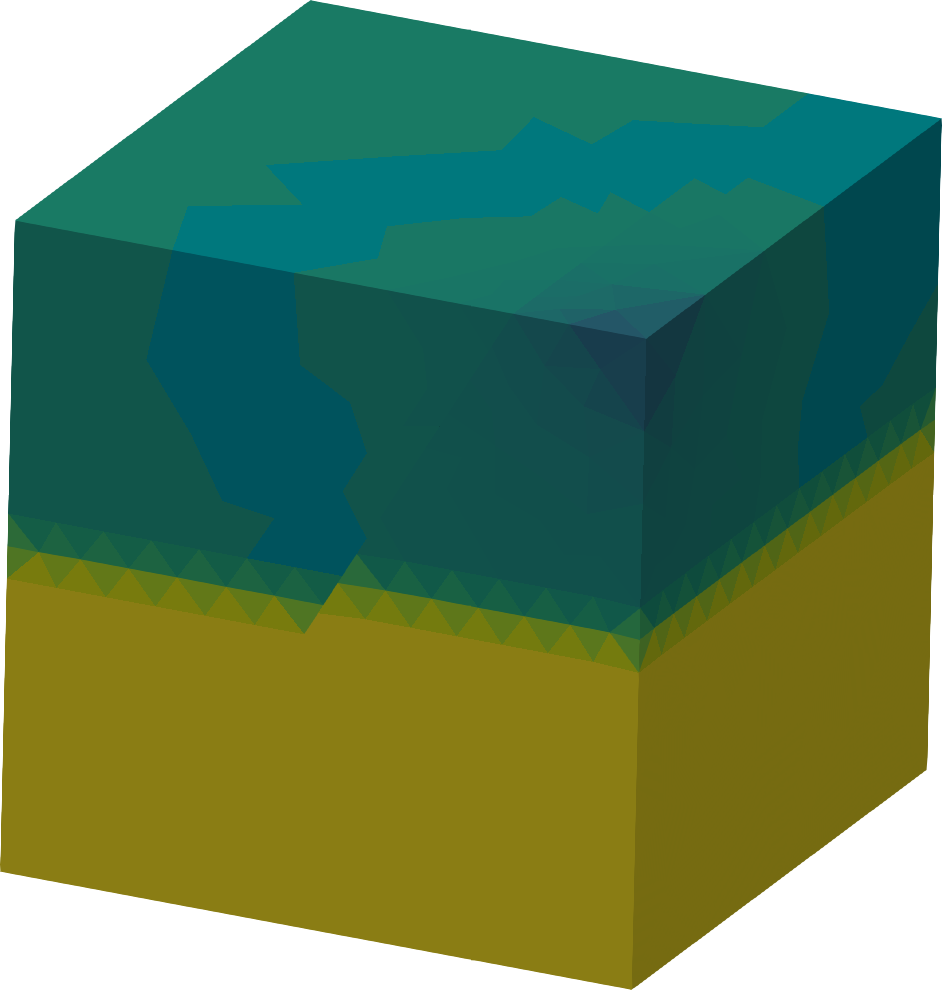}}
    \\
    \subfloat[block-POD]{\includegraphics[width=0.35\textwidth]{ 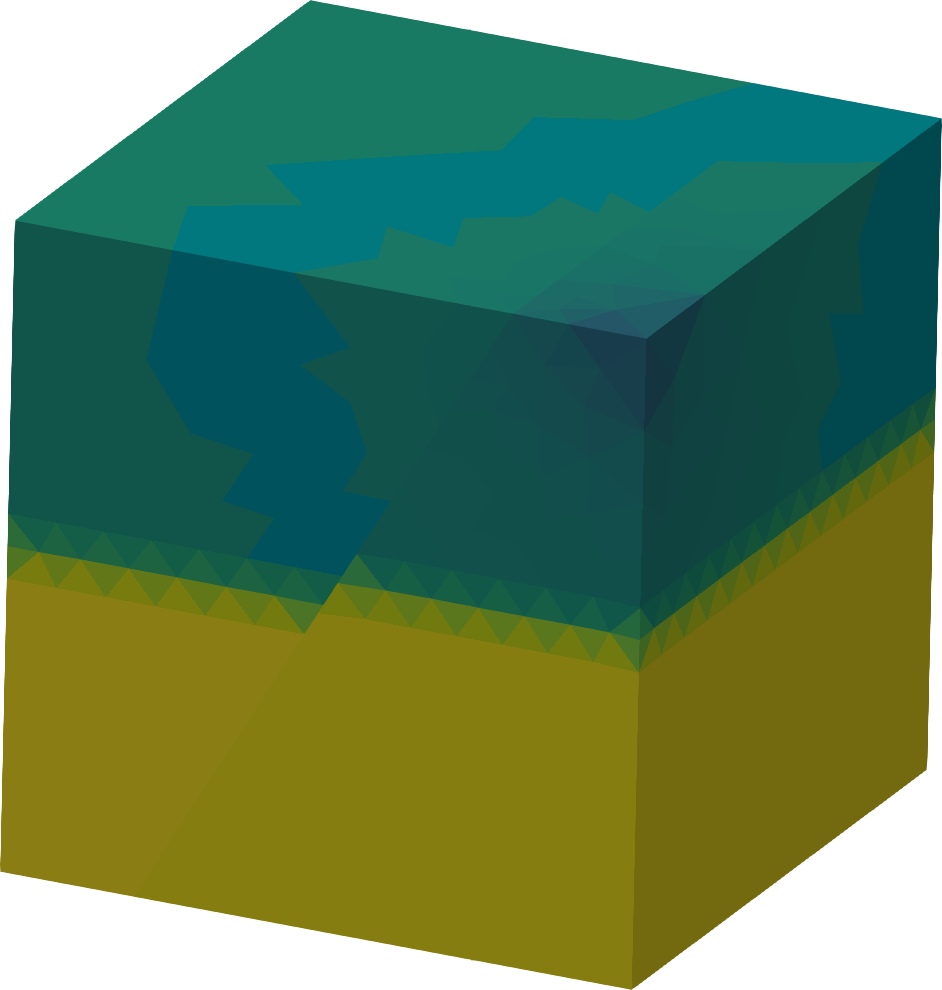}}
    \hspace*{0.1cm}
    \subfloat[DL-ROM]{\includegraphics[width=0.35\textwidth]{ 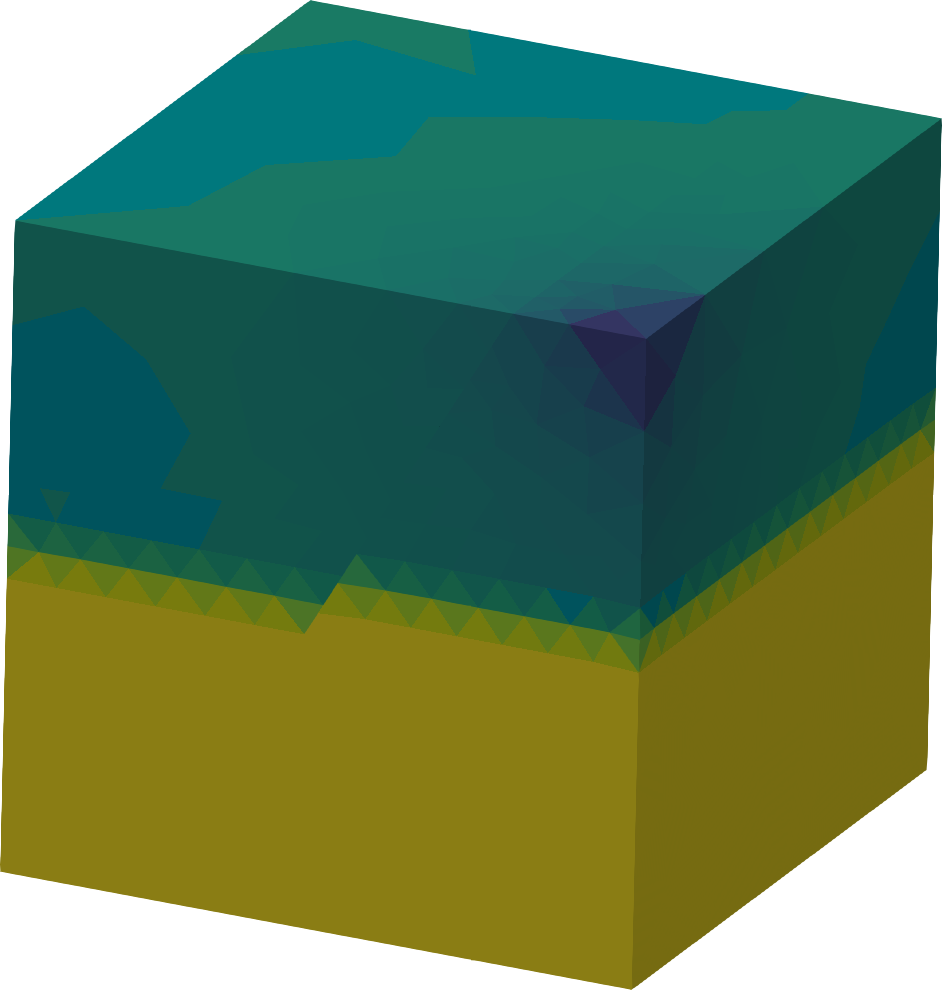}}
    \vspace{0.3cm}
    \includegraphics[width=0.35\textwidth]{ 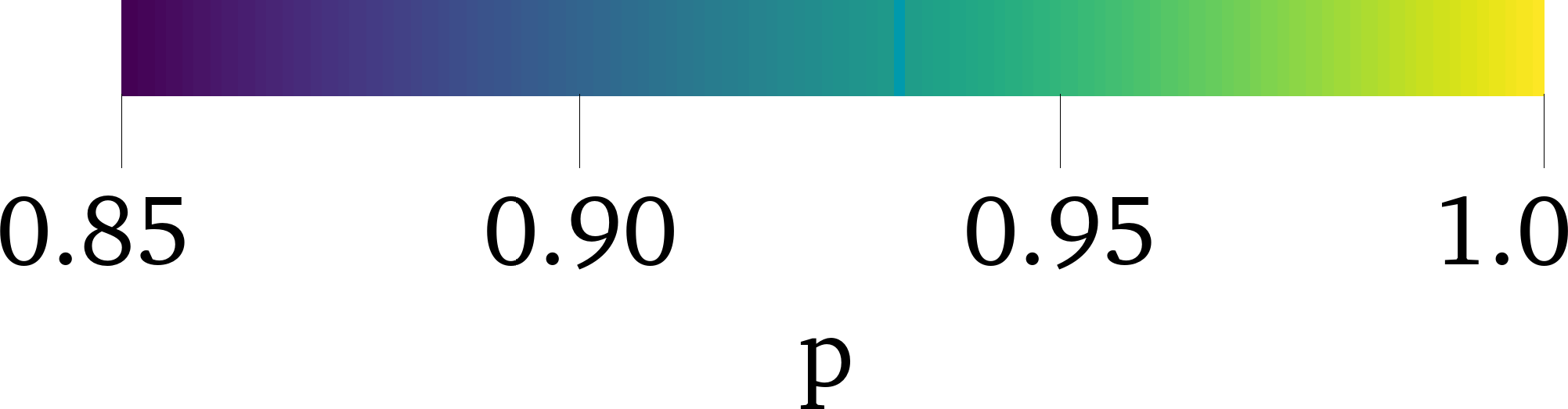}
    \caption{Case study 2. Reconstructed solutions for a specific value of the parameter $\mu$.}
    \label{fig:test_2_single_snap}
\end{figure}

\begin{figure}
    \centering
    \includegraphics[width=7cm]{ 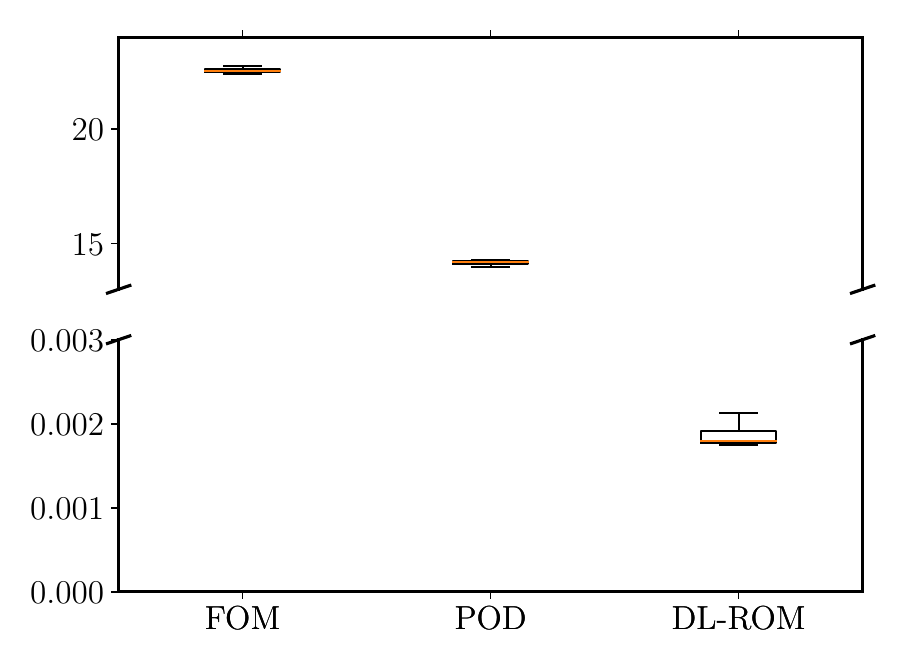}
    \caption{Case study 2.  Evaluation times in seconds of the full-order model, reduced model obtained through POD, and reduced model obtained through the DL-ROM.}
    \label{fig:test_2_online_time}
\end{figure}

The specific architectures are summarized in Tab.~\ref{tab:test_2_dense}. We can infer that the compression rate is equal to $\frac{n}{N} = 0.06\%$, which is ten times higher than the POD.
The hyperparameters of the training are the same as in Case 1, except for the number of epochs that is higher, 6000 instead of 4000 (about 3 hours).
The errors in the test data set are: $e_{max} = 4.6\%$, $e_{min} = 0.16\%$, $e_{ave} = 1.02\% $. They are similar to the values of case 1, suggesting that the DL-ROM works equally in a 2D and 3D scenario.
Fig.~\ref{fig:test_2_single_snap} shows a snapshot of the test data set. As in \cref{sec:test_1}, we can see that, qualitatively, all the solutions resemble the FOM one at each point of the domain.  The online time is reported in Fig.~\ref{fig:test_2_online_time}, despite the large increase in degrees of freedom of the problem, the online time has only doubled with respect to case 1. The consequent time saving is equal to $99.97\%$.
\begin{table}[]
\begin{tabular}{llll}
                                     & layer        & \# tw       & \# tw tot    \\
                                     \hline \\
\multirow{3}{*}{Encoder}             & fc, 8706, 1708, PReLU & \multirow{2}{*}{\numprint{14880103}}  & \multirow{8}{*}{\numprint{29 770 013}} \\
                                     & fc, 1708, 5, PReLU    &                            &                             \\
                                     &                       &                         &                           \\
\multirow{3}{*}{Decoder}             & fc, 5, 1708, PReLU    & \multirow{2}{*}{\numprint{14888803}} &                             \\
                                     & fc, 1708, 8706, PReLU &                            &                             \\
                                     &                       &                         &                           \\
\multirow{2}{*}{Reduced map network} & fc, 5, 100, PReLU     & \multirow{2}{*}{1107}      &                             \\
                                     & fc, 100, 5, PReLU     &                            & \\
                                     \hline
\end{tabular}
\caption{Case 2. Architecture of the networks that make up the reduced model. ``fc" stands for fully connected layer, tw is the number of trainable weights.}
\label{tab:test_2_dense}
\end{table}
%
%
%

\subsection{Case 3 - setup}\label{sec:test_3}
\begin{figure}[h]
    \centering
    \includegraphics[width=0.4\textwidth]{ 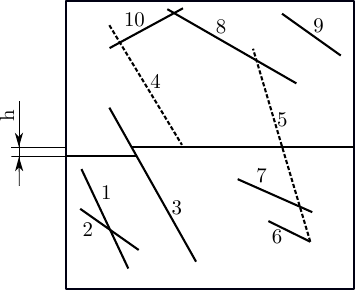}
    \caption{Case study 3. The domain contains 9 fractures and one fault (number 3). Fractures 4 and 5 are blocking fractures, while all others have high permeability.}
    \label{fig:test_3_domain}
\end{figure}

\begin{figure}[h]
    \centering
    \subfloat[FOM]{\includegraphics[width=0.3\textwidth]{ 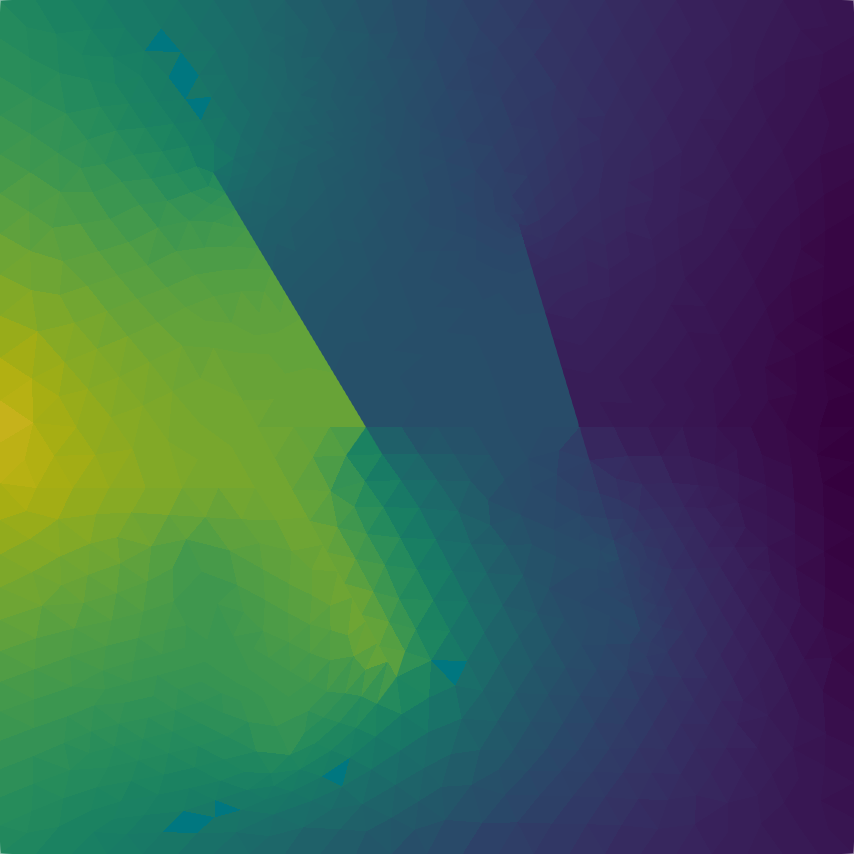}}
    \hspace*{0.01\textwidth}
    \subfloat[FOM]{\includegraphics[width=0.3\textwidth]{ 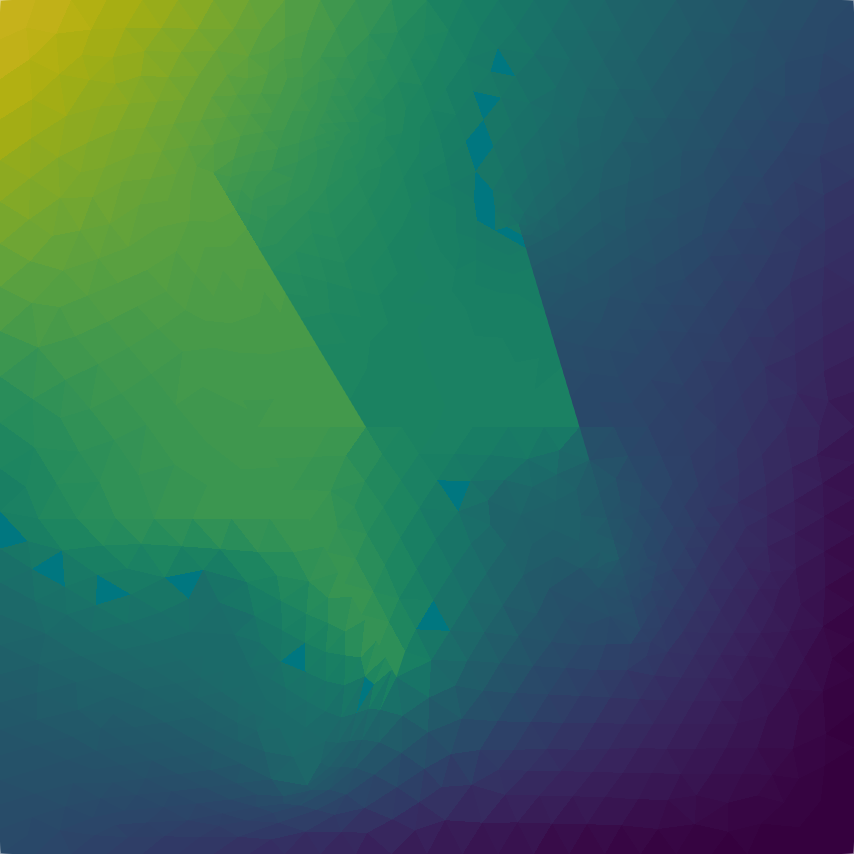}}
    \hspace*{0.01\textwidth}
    \subfloat[FOM]{\includegraphics[width=0.3\textwidth]{ 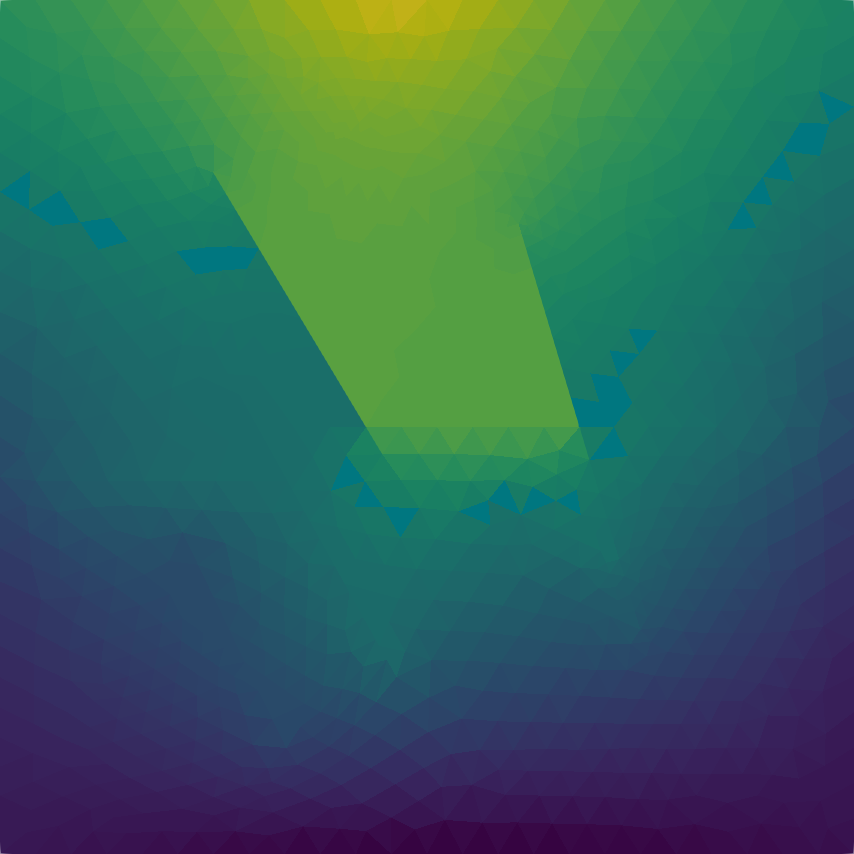}}
    \vspace{0.2cm}
    \includegraphics[width=0.3\textwidth]{ 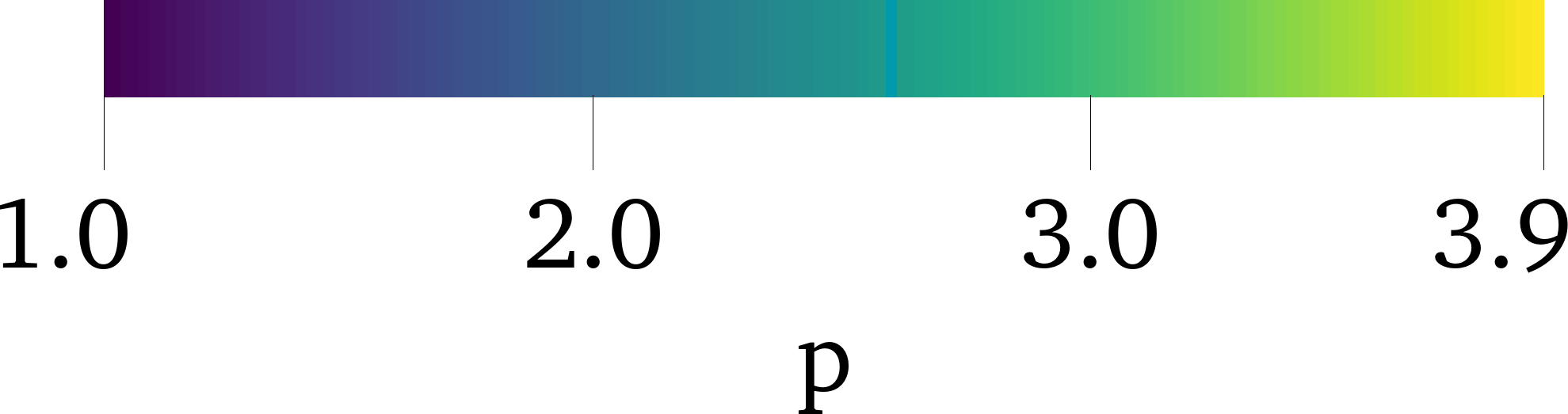}
    \caption{Case study 3. Boundary conditions are applied such that the mean pressure gradient is (a) horizontal (b)  slanted (c) vertical.}
    \label{fig:test_3_single_snap_fom}
\end{figure}

The third test case resembles the complex fracture network benchmark in \cite{Berre_Boon2021}. The 2D square domain includes intersecting permeable and impermeable fractures, as shown in Fig.~\ref{fig:test_3_domain}, whose detailed description can be found in the reference. Compared to the reference, we further increase the complexity by setting two horizontal layers whose rock permeability is equal to $10^{-2}$ at the bottom and $10^2$ at the top.
A Dirichlet boundary condition for the pressure is applied to all the boundaries such that an average flow is generated from one side to the opposite. After a manual sensitivity analysis, we set the two most relevant features as uncertain parameters: the first controls the boundary condition creating high variation among the snapshots, see Fig.~\ref{fig:test_3_single_snap_fom}, while the second controls the geometry, making the problem non-affine. A sinusoidal variation of the pressure at the boundary, $p_b$, is applied:
\begin{equation*}
p_b(\omega) = p_1 \left(1-\sin\left(\frac{\omega - \omega_0}{2}\right)\right) + p_2 \left( \sin\left(\frac{\omega - \omega_0}{2}\right) \right),
\end{equation*}
where $p_1$ and $p_2$ are the maximum and minimum values, $\omega = \arctan(y/x)$, $x$ and $y$ are the coordinate of the boundary points.
The first parameter is the reference angle $\omega_0$ which can be subject to a variation of $90\degree$,  Fig.~\ref{fig:test_3_single_snap_fom} shows snapshots for three different values of $\omega_0$. The second parameter controls the position of the horizon on the left of fault number 3, displacing rigidly also fractures number 1 and 2.
The high variation of the boundary conditions and the presence of blocking fractures imply strong nonlinearities in the solution manifold that may raise issues in reducing the order of the model.
Generally speaking, the results of this case show a larger gap between the performance of the two reduced-order model methods in favor of the DL-ROM one.

\subsubsection{Case 3 - POD}
\begin{figure}[h]
    \centering
    \includegraphics[width=5cm]{ 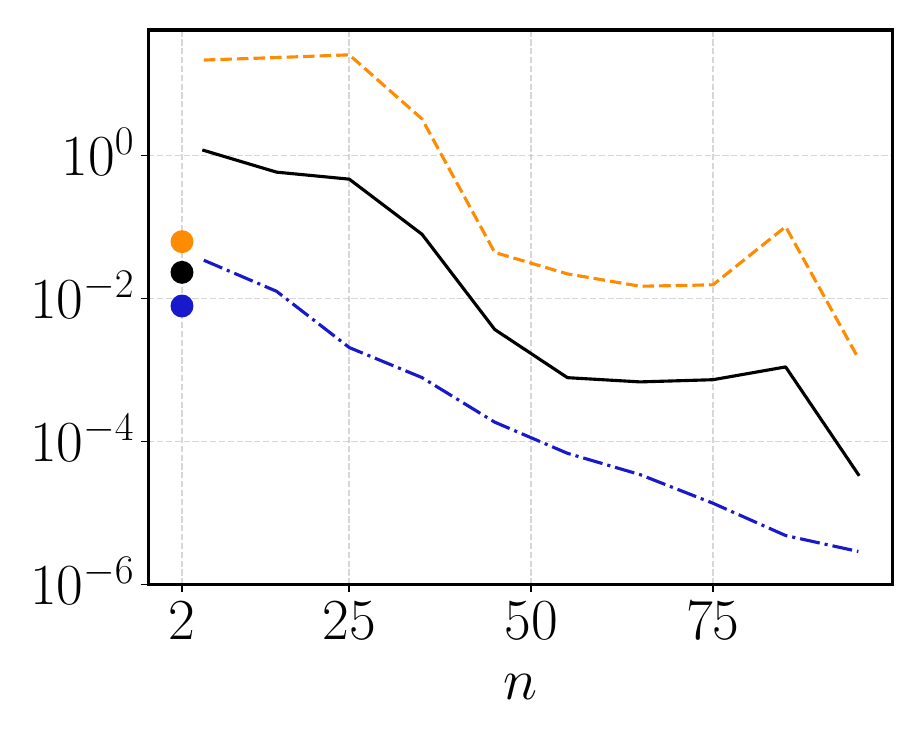}
    \vspace{0.2cm}
    \includegraphics[width=10cm]{ test_1_legend.pdf}
    \caption{Case study 3. $e_{max}$, $e_{min}$, $e_{ave}$, respectively, maximum, minimum and average error versus the size of the reduced space, $n$, of the reduced model obtained with the standard POD. Points represent the value errors for the DL-ROM, described in \Cref{sec:test_2_dlrom}.}
    \label{fig:test_3_err_all}
\end{figure}

For this test, only the standard POD is investigated. 
The error versus the size of the reduced space (Fig.~\ref{fig:test_3_err_all}) shows that several modes are required to reach a low error.
The online time is depicted in Fig.~\ref{fig:test_3_online_time} where a number of $n = 44$ modes is used. The strong nonlinearities of the problem manifest in a compression rate higher than in the previous cases $\frac{n}{N} = 1.85\%$. The five-number summary of the online time is shown in Fig.~\ref{fig:test_3_online_time}.
%
%

\subsubsection{Case 3 - DL-ROM}
\begin{figure}[h]
    \centering
    \includegraphics[width=7cm]{ 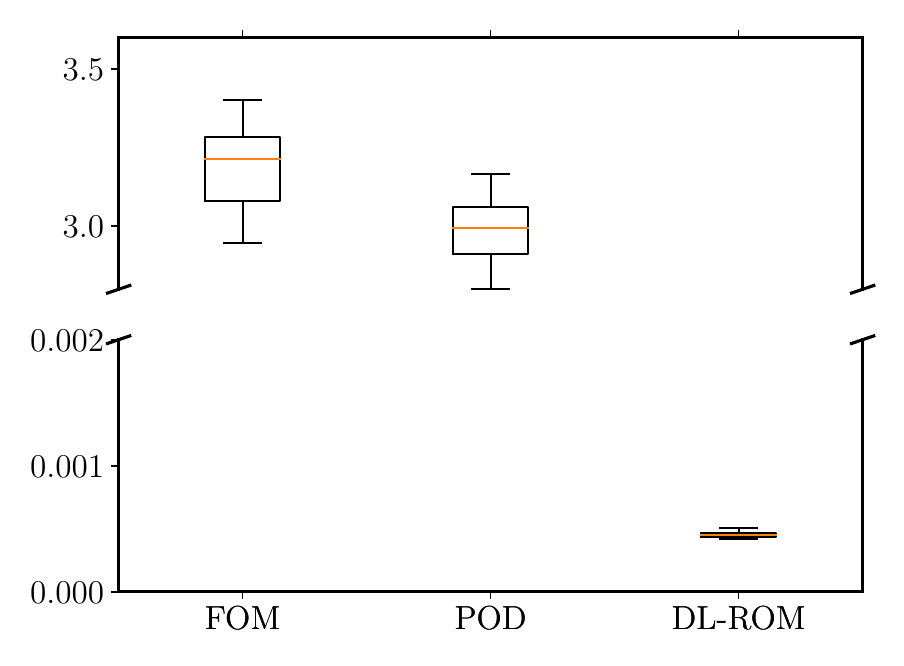}
    \caption{Case study 3. Evaluation times in seconds of the full order model, reduced model obtained through POD, and reduced model obtained through the DL-ROM.}
    \label{fig:test_3_online_time}
\end{figure}

The architectures are described in Tab.~\ref{tab:test_3_dense}, in particular, the reduced space size is equal to $n = 2$, so the compression rate results in $\frac{n}{N} = 0.08\%$.
The training requires 5000 epochs (about 74 minutes) to reach a steady low value of the loss function. The errors are fairly low: $e_{max} = 6.2\%$, $e_{min} = 0.78\%$, $e_{ave} = 2.3\%$ (Fig.~\ref{fig:test_3_err_all}), although the ratio between the minimum and maximum values is higher than the one in case 1 and case 2.
The online time can be inferred from Fig.~\ref{fig:test_3_online_time}, the time savings are equal to $99.86\%$.
\begin{table}[]
\begin{tabular}{llll}
                                     & layer            & tw     & tw tot \\
                                     \hline \\
\multirow{3}{*}{Encoder}             & fc , 2378, 404, PReLU & \multirow{2}{*}{\numprint{961928}}  & \multirow{8}{*}{\numprint{1926735}} \\
                                     & fc, 404, 2, PReLU         &                          &                          \\
                                     &                       &                         &                           \\                                     
\multirow{3}{*}{Decoder}             & fc, 2, 404, PReLU         & \multirow{2}{*}{\numprint{964 303}} &                          \\
                                     & fc, 404, 2378, PReLU  &                          &                          \\
                                     &                       &                         &                           \\                                     
\multirow{2}{*}{Reduced map network} & fc, 2, 100, PReLU         & \multirow{2}{*}{504}     &                          \\
                                     & fc, 100, 2, PReLU         &                          &                          \\
                                     \hline 
\end{tabular}
\caption{Case study 3. Architecture of the networks that make up the reduced model. ``fc" stands for fully connected layer, tw is the number of trainable weights.}
\label{tab:test_3_dense}
\end{table}
%
%
%
\begin{figure}[h]
    \centering
    \subfloat[FOM]{\includegraphics[width=0.3\textwidth]{ test_3_sol_fom_969.png}}
    \hspace*{0.01\textwidth}
    \subfloat[POD]{\includegraphics[width=0.3\textwidth]{ 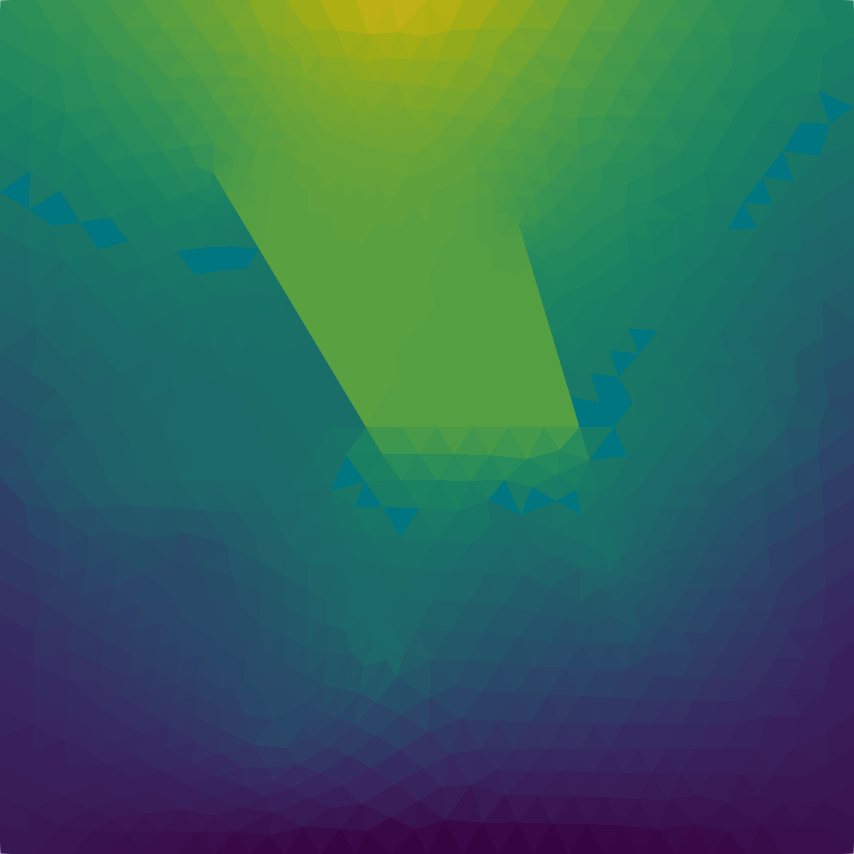}}
    \hspace*{0.01\textwidth}
    \subfloat[DL-ROM]{\includegraphics[width=0.3\textwidth]{ 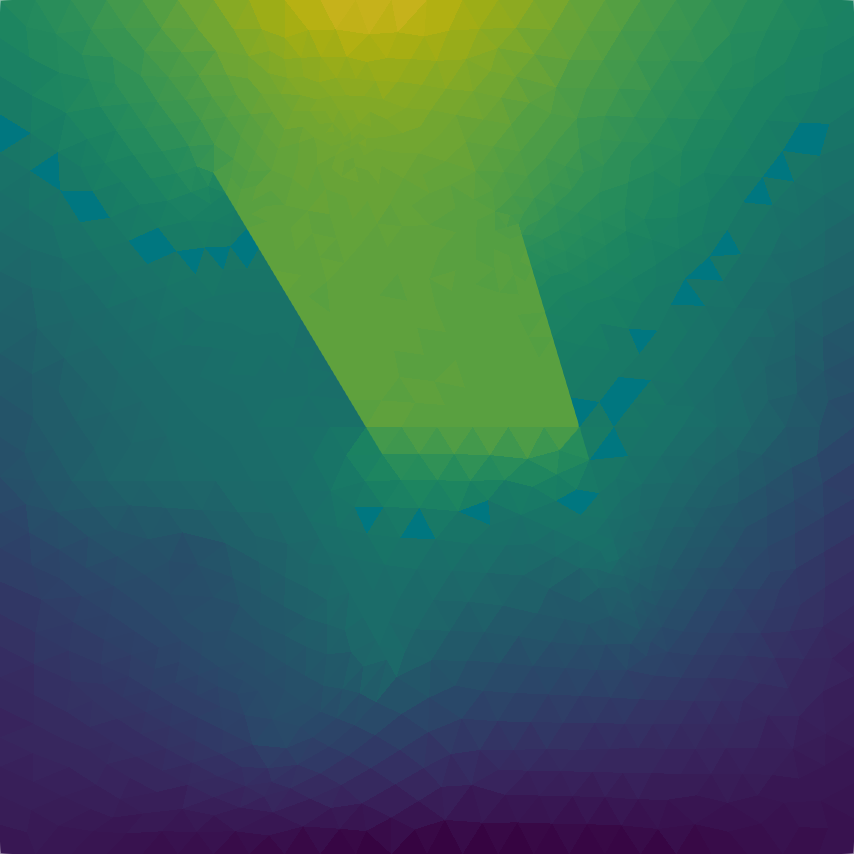}}
    \vspace{0.2cm}
    \includegraphics[width=0.3\textwidth]{ test_3_pressure_label.png}
    \caption{Case study 3. Reconstructed solutions for a specific value of the parameter $\mu$.}
    \label{fig:test_3_single_snap}
\end{figure}
%

%

\section{Multi-query application}\label{sec:multi-query}
In this section, we show possible practical applications of the reduced order model techniques. We first used a Monte Carlo strategy to sample the data and conduct a sensitivity analysis, comparing the results of different reduced-order models and a full-order model. Afterwards, an inverse problem is solved exclusively with DL-ROM.
%

\subsection{Sensitivity analysis}\label{sec:sensitivity}
We use a Monte Carlo technique to perform a sensitivity analysis on a problem based on test case 1 (see \Cref{sec:test_1}). We chose this test case because its small size allowed us to run the Monte Carlo analysis even with the full-order model and thus compare the results with those obtained with the reduced models.
We rely on the Chaospy library, a numerical toolbox to perform uncertainty quantification \cite{Feinberg2015}.
Our goal is to focus on a quantity of interest (q.o.i) represented by the pressure jump between injection and production, see Fig.\ref{fig:test_1_domain}, called $\Delta p$, and determine its mean value, $\overline{\Delta p}$, standard deviation $\Tilde{\sigma}$, and the first-order sensitivity index (Sobol index) \cite{Saltelli2004}, $\Tilde{s}_{1,i}$, $i=1,\ldots,5$, related to each entry, $\mu_i$, of the parameter vector $\mu$.

For the purpose of this example, we prepare a training data set of 200 snapshots and a validation data set of 20 samples. All data in the training data set are used to create the snapshot matrix to which the SVD decomposition is applied. Then, the left singular vector matrix, $U$, see \Cref{sec:pod}, is truncated to $n = 13$ basis functions because it leads to a low enough reconstruction error, that is, $e_{ave} = 5\%$. The neural networks are trained until the loss function reaches a low enough value to have a satisfactory reconstruction error equal to $e_{ave} = 5.3\%$. Since the q.o.i. is known a priori, it is possible to add a term in the loss function related to the q.o.i., thus increasing the accuracy of the reduced model for only what is needed. For instance, in this case, we add the term $\ell_3$ to the loss function defined in \eqref{eq:loss_function}:
\begin{align*}
    \mathscr{L} &= \alpha \ell_1 + \beta \ell_2 + \gamma \ell_3  = \\
    &= \frac{\alpha}{N} \| u_N(\mu_i) - \Psi(\Psi'(u_N(\mu_i))) \|_2^2 + \frac{\beta}{n} \| \Psi'(u_N(\mu_i)) - \varphi(\mu_i) \|_2^2 + \gamma ( \Delta p_{rom} - \Delta p_{fom}  )^2 
\end{align*}
where $\gamma$ is a user-defined coefficient, $\Delta p_{rom}$ and $\Delta p_{fom}$ are, respectively, the jump of pressure between the injection and the production computed with the reduced order model and the full order model.

\begin{figure}[h]
    \centering
    \includegraphics[width=0.75\textwidth]{ 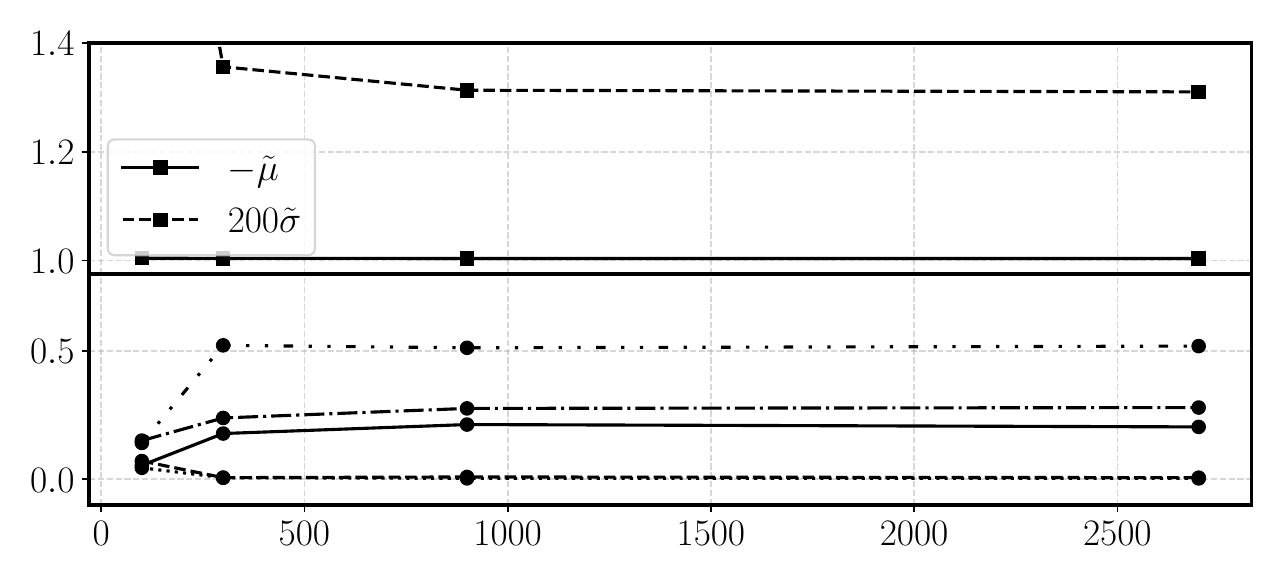}
    \includegraphics[width=0.75\textwidth]{ 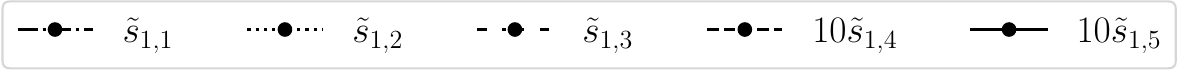}
    \caption{Convergence of relevant statistics. The values reach a stable value with a set of 900 samples or more.}
    \label{fig:statistics_convergence}
\end{figure}

In order to reach stable values of the statistics, having a negligible error due to the finite sampling, 
we generate a number of 900 instances from a uniform distribution of the uncertain parameters to compute the desired statistics. This value is inferred from convergence analysis, made computationally cheap by the neural network approach, shown in Fig.~\ref{fig:statistics_convergence}.

\begin{figure}[h]
    \centering
    \includegraphics[width=7cm]{ 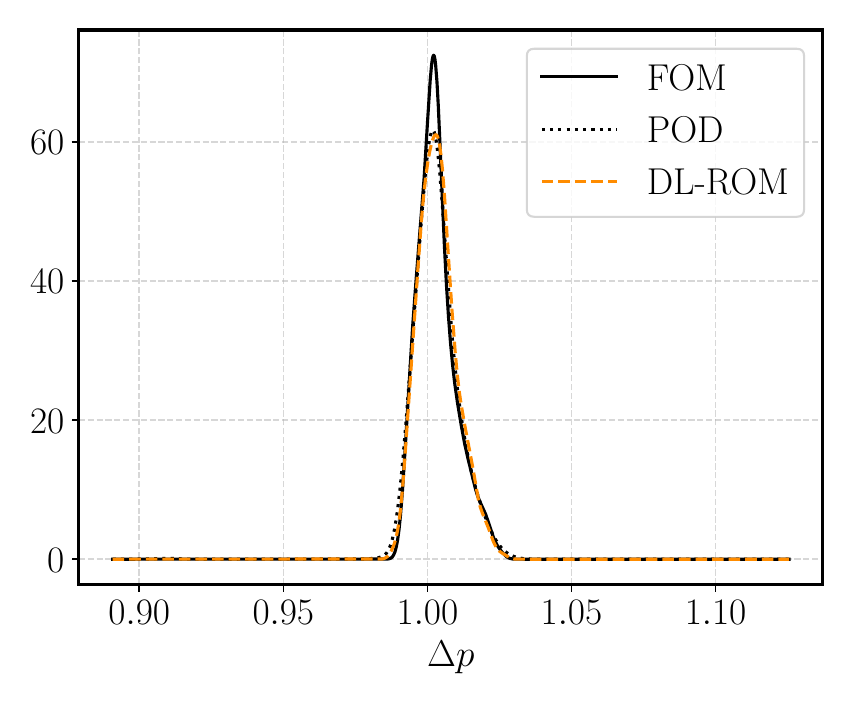}
    \caption{Probability density function of the quantity of interest, $\Delta p$. Both POD and DL-ROM well reproduce the shape.}
    \label{fig:pdf_qoi}
\end{figure}

Fig.~\ref{fig:pdf_qoi} shows the probability density function (pdf) of the q.o.i. estimated using a Gaussian kernel estimator \cite{Scott1992}. Its profile is well reproduced by both reduced models. 

We list in Tab.~\ref{tab:statistics} the summary statistics of our interests. We notice a good agreement between the values; most significant discrepancies occur for small quantities, such as $\Tilde{s}_{1,1}$ and $\Tilde{s}_{1,4}$, whose values are difficult to capture accurately. From the Sobol index, we can see that most of the variance in the pressure jump is due to the permeability of layer 1, the bottom layer, and layer 3, the caprock, since their permeabilities are low. The small displacement along the fault is less relevant, and the high-permeability upper layer and fault permeability play a minor role in this context. 

We now want to estimate the total time required to assess this multi-query application. A single numerical solution of the full-order problem takes on average $0.58 \ s$  per core, see Fig.~\ref{fig:test_1_online_time}, using a 4 core CPU, we get a total time for the Monte Carlo evaluation of about $131 \ s$. For the POD we need to consider the offline time to create the reduced model and the time to generate the required samples with the reduced model. Given the size of the training data set, we see that data generation takes $200\times0.58/4 = 29 \ s$. The creation of the $\Phi$ matrix requires a negligible time compared to the other operations. The online time is $0.57 \ s$, which is slightly less than the full-order model time because the latter is already small and most of the time is spent deforming the mesh and reassembling the full-order matrix $A$. The POD sampling time results to be equal to $900\times0.57/4 = 130 \ s$, therefore, the total time to evaluate the sensitivity analysis is $159 \ s$ plus an additional time for additional routines that are common for all methods.
Training in the neural network takes $400$ epochs in $67 \ s$. The overall sampling time is less than one second, so it is negligible. The total time, which includes the training and validation data sets, is $96 \ s$. 
This problem is very computationally cheap, so a reduced-order model technique does not show great computational time improvements, indeed the POD takes even longer than the full-order model. We recall that we chose this problem for the possibility of running the Monte Carlo analysis with the full-order model, hence to compare the accuracy of the results.
Moreover, we repeat that those are just indicative numbers; the computational time depends on the relation between the algorithm, code, and hardware since different parts of the algorithm run better on different hardware because of different code parallelization. Depending on the user's hardware availability, the outcomes may vary, but the neural network approach has an online time order of magnitude faster than the other methods, so for a problem large enough it will be convenient.
\begin{table}
{\footnotesize
\begin{tabular}{lccccccc}
       & $\overline{\Delta p}$ & $\Tilde{\sigma}\times10^3$   & $\Tilde{s}_{1,1}\times10^{1}$ & $\Tilde{s}_{1,2}$ & $\Tilde{s}_{1,3}\times10^{1}$ & $\Tilde{s}_{1,4}\times10^{3}$ & $\Tilde{s}_{1,5}\times10^{2}$ \\
       \hline
FOM    & 1.0033           & $6.61$  & $3.11$ & $7.04\times10^{-8}$ & $5.86$ & $0.32$ & $0.85$ \\
POD    & 1.0028           & $7.79$  & $1.28$ & $8.16\times10^{-3}$ & $3.62$ & $9.37$ & $2.45$ \\
DL-ROM & 1.0036           & $6.57$  & $2.76$ & $3.73\times10^{-3}$ & $5.12$ & $0.92$ & $2.13$ \\
\hline
\end{tabular}
}
\caption{Mean value, $\overline{\Delta p}$, deviation, $\Tilde{\sigma}$, and first order sensitivity index, $\Tilde{s}_{1,i}$, of the quantity of interest obtained with data generated by the FOM, POD, and DL-ROM.}
\label{tab:statistics}
\end{table}

\subsection{Inverse problem}

The goal of this application is to determine the parameters such that a q.o.i. is equal to a desired value. We take Case 2 and set the pressure difference between injection and production equal to the desired value of $\Delta p_{d} = 0.188$ derived from the following set of parameters: $K_1 = 0.1$, $K_2 = 150$, $K_3 = 1\times10^{-4}$, $K_4 =9.5\times10^{-4}$, $h = 0.09$.
Assuming that some parameters, $K_2$ and $K_4$, are fairly known while some others are more uncertain, the ranges where we seek a solution are: $K_1 \in [10^{-4}, 1]$, $K_2  \in [10^{2}, 2\times10^{2}]$, $K_3 \in [10^{-6}, 10^{-4}]$, $K_4 \in [9\times10^{-4}, 10^{-3}]$, and $h \in [0.01, 0.1]$.

In this example, we show only the application of DL-ROM. Similarly to the previous application, neural networks are trained on a small dataset made up of $200$ snapshots because high precision is not required for the purpose of the current application. The training data set contains solutions sampled from the ranges defined in \ref{sec:test_1} that are stricter than those considered in this application, so we indirectly exploit the extrapolation capacity of the reduced model.

We cast the inverse problem as: $\min_{\mu} F(\mu)$, where $F(\mu) = (\Delta p(\mu) - \Delta p_{d})^2$.
We want to show that DL-ROM allows us to use heuristic optimization algorithms, which usually require a large number of evaluations of the objective function. We select the differential evolution algorithm \cite{Storn1997} implemented in Scipy \cite{Virtanen2020} with the default setting, except for the convergence tolerance: $tol = 0.001$ and $atol = 10^{-10}$.
After $300$ iterations, with a total number of objective function evaluations of $22581$, the algorithm satisfies the convergence criterion. The optimal solution gives $K_1 = 0.251$, $K_2= 153$, $K_3=9.1\times10^{4}$, $K_4=9.4\times10^{-4}$, and $h = 0.08$. We observe that the optimal values of the parameters are similar to the exact ones except for $K_1$ and a small error affects also the value of $h$.

The $\Delta p$ obtained with the reduced model is equal to $0.188$ as requested, while taking the optimal solution, recomputing the results of $\Delta p$ with the results of the full-order model equal to $0.196$, close enough to the desired value.

\begin{figure}[h]
    \centering
    \includegraphics[width=7cm]{ 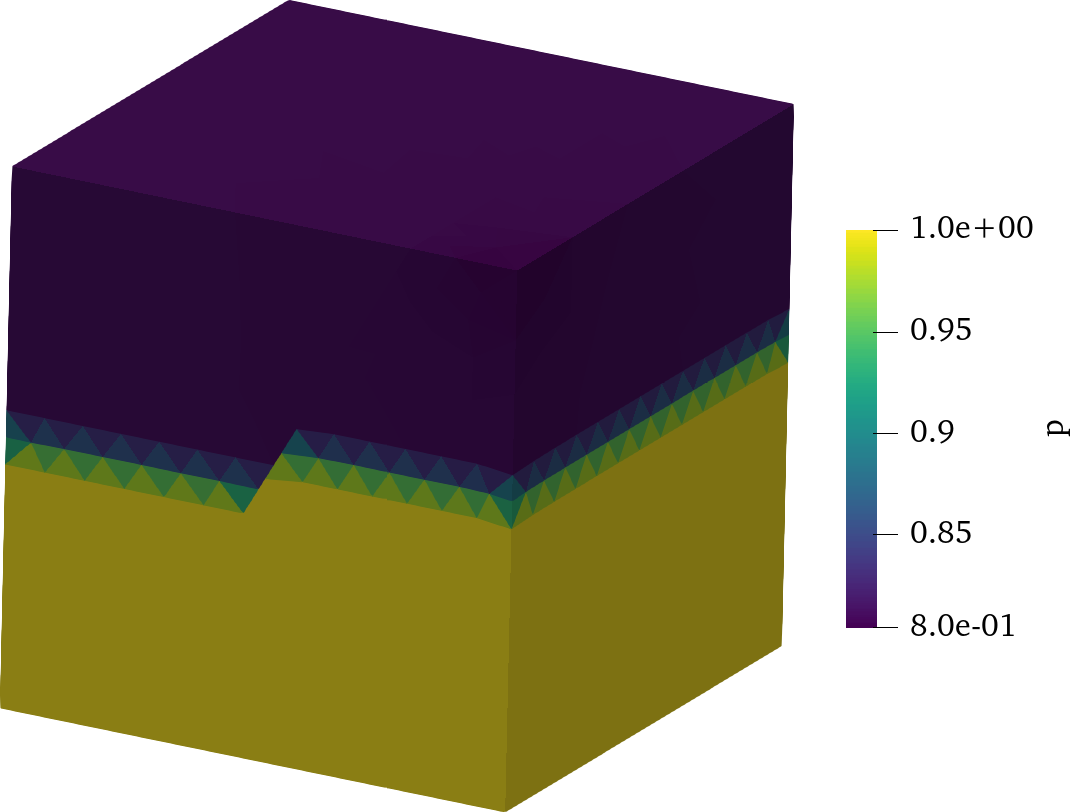}
    \caption{Full order model pressure field obtained from the optimal solution.}
    \label{fig:optimal_sol}
\end{figure}

See Fig.~\ref{fig:optimal_sol} for the full-order model pressure field obtained with the optimal solution.

Due to the high number of objective function evaluations, we have considerable time savings: the offline phase takes about $75$ min for snapshot generation and $4$ min for neural network training, while $22581$ runs would require $142$ hours with the FOM but only $136$ s with the use of the DL-ROM. Therefore, the total computational time is approximately $142$ hours without the use of reduced models, and $1$ h $18$ min with the DL-ROM.

%

\section{Conclusion}

We focussed on reduced-order modeling techniques applied to the problem of a single-phase flow in rigid porous media with an arbitrary number of fractures and faults. The mixed-dimensional framework lets us efficiently deal with geometrical discontinuities (fractures or faults) in a geometry-deforming setup. It was possible to generate the data (solution of the full-order model) and create the reduced-order model without any further devices. We considered uncertainties in the parameters with respect to the physical properties of the rock and the geometry, whereas uncertainties in the fluid properties represent interesting possible future developments. 

According to the results shown by our tests, the DL-ROM takes slightly longer compared to classic methods as the POD to generate the reduced model because of the training of the neural networks, but the online time is extremely low, which makes this approach promising, especially when the number of queries required is large.

Similarly to block-POD, the efficiency of the DL-ROM could be improved by segregating variables and severing the links between disparate physical variables, i.e., $p$, $p_\gamma$, $\lambda$. This would decrease the total trainable weights, thereby expediting the training phase. Although we have conducted an initial review of this approach, it remains under development and is the subject of future research.

A disadvantage of the neural network approach is the large number of hyperparameters, such as the size of the training dataset, architecture of the neural networks, parameters of the optimization algorithm, etc.,  for which some user experience is needed since they affect the accuracy of the reduced model. We showed that the ROM strategies, and in particular the DL-ROM, lead to an advantage in terms of faster analysis with satisfying accuracy, so further investigations will be undertaken on the line of multifidelity ROM, and study of applications of ROM strategies to time-dependent problems simulating more realistic and complex physics. For instance, this could involve studying scenarios such as two-phase flows in fractured porous media, which presents, in addition to the geometric discontinuities, the challenge due to the sharp fronts resulting from the hyperbolic nature of the problem. This highly nonlinear scenario may pose strong difficulties for methods based on a linear map, $\mathcal{M}$, while promoting nonlinear methods such as the DL-ROM.

%

\section{Acknowledgments}

The research has been supported by the Italian Ministry of Universities and Research (MUR) under the project ``PON Ricerca e Innovazione 2014-2020" and was carried out in collaboration with Eni S.p.A.
The last three authors also gratefully acknowledge the support of the "Dipartimento di Eccellenza 2023-2027".
All authors warmly thank Andrea Manzoni, Stefano Micheletti, and Nicola Rares Franco
for many insightful discussions.

%

%

\small

\bibliography{bibliography, biblio}

%

\appendix
\section{Nomenclature}\label{appendix:nomenclature}

\begin{longtable} {p{.30\textwidth} p{.70\textwidth}}
    $a$ & Trainable parameter of PReLU \\
    $A$ & Matrix describing the discrete equation \\
    $b$ & Discrete right-hand side \\
    $C, C_d, C_{df}, C_s, C_{sf}$ & Control points sets \\
    $\overline{C}_s$ & Set of surfaced where sliding conditions are applied \\
    $d(x,y)$ & Euclidean distance between $x$ and $y$ \\
    $D$ & Spatial dimension \\
    $e$ & Number of parameters \\
    $e_{min}, e_{max}, e_{ave}$ & Minimum, maximum, averaged relative errors between the FOM and ROM solution \\
    $f$ $\sib{1\per\second}$ & Scalar source or sink term \\
    $f_\gamma$ $\sib{1\per\second}$ & Scalar source or sink term in the fault \\
    $G$ & Matrix of displacement constraint \\
    $g$ & Radial basis function dependent on the distance, $d$, of two points, $x_1, x_2$ \\
    $g^*$ & Radial basis function dependent on $x_1, x_2$ \\
    $g^\dag$ & Modified radial basis function \\
    $H_k$ & Matrix of no tangential contribution constraint applied to surface $k$ \\
    $\mathcal{I}$ & influence function \\
    $l$ & Number of control points \\
    $N$ & Number of degrees of freedom of the full order problem \\
    $n$ & Number of degrees of freedom of the reduced problem \\
    $n_{br}$ & number of intersecting branches \\
    $n_{\min}$ & Minimal latent dimension \\
    $p$ $\sib{\pascal}$ & Pressure \\
    $\overline{p}$ $\sib{\pascal}$ & Pressure on boundaries \\ 
    $p_\gamma$ $\sib{Pa}$ & Pressure in the fault \\
    $p_\iota$ $\sib{Pa}$ & Pressure at the intersection \\
    $\overline{\Delta p}$ $\sib{Pa}$ & Mean value of $\Delta p$ \\
    $q$ $\sib{\meter\per\second}$ & Darcy velocity \\
    $\overline{q}$ $\sib{\meter\per\second}$ & Darcy velocity on boundaries \\
    $K$ $\sib{\meter^3\second\per\kilo\gram}$ & Intrinsic permeability scaled by the dynamic viscosity \\
    $K_n$ $\sib{\meter^3 \second \per \kilo\gram}$ & Normal fault permeability \\
    $K_\tau$ $\sib{\meter^3 \second \per \kilo\gram}$ & In-plane fault permeability \\
    $K_\iota$ $\sib{\meter^3 \second \per \kilo\gram}$ & Representative permeability at intersection \\
    $r$ & Map from $\gamma$ to $\partial_{in}\Omega$ \\
    $S$ & Snapshot matrix \\ 
    $s$ $\sib{\meter}$ & Displacement \\
    $\overline{s}$ $\sib{\meter}$ & Known displacement \\
    $\Tilde{s}$ & First order sensitivity index \\
    $t, b$ & Non-parallel tangent unit vectors of sliding surface \\
    $u_N$ & Full order model solution \\
    $\Tilde{u}_N$ & Reconstructed solution \\
    $u_n$ & Reduced order model solution \\
    $U$ & Left singular vector matrix \\
    $U_{tr}$ & Left singular vector matrix truncated \\
    $V$ & Right singular vector matrix \\
    $z$ & Unknown of mesh deformation linear system \\
    $\mathscr{L}$ & Loss function \\
    $\mathcal{M}$ & Map from full order model space to reduced space\\
    $\mathcal{S}$ & Solution manifold \\
    $\mathcal{V}_n$ & Reduced problem solution space \\
    $\mathcal{V}_N$ & Full order model solution space \\
    $\alpha, \beta$ & User-defined loss function weights \\
    $\beta$ & Side function \\
    $\gamma$ & Fault domain \\
    $\partial_p \gamma$ & Boundary of $\gamma$ where Dirichlet boundary condition for the pressure is applied \\
    $\partial_q \gamma$ & Boundary of $\gamma$ where Neumann boundary condition is applied \\
    $\partial_{ex}\gamma$ & Boundary of $\gamma$ in contact with $\partial\Omega$ \\
    $\partial_{in}\gamma$ & Boundary of $\gamma$ not in contact with $\partial\Omega$ \\
    $\gamma^+ (\gamma^-)$ & Additional interfaces between the matrix domain, $\Omega$ and fault domain, $\gamma$ \\
    $\delta_n$ & Non-linear counterpart of Kolmogorov $n$-width. \\
    $\epsilon$ $\sib{\meter}$ & Fault aperture \\
    $\zeta$ & Unknown coefficients of the linear combination of radial functions \\
    $\eta$ & Exponent defining the permeability \\
    $\lambda^+ (\lambda^-)$ $\sib{\meter\per\second}$ & Volumetric fluid flux exchanged between subdomains \\
    $\lambda_\gamma$ $\sib{\meter\per\second}$ & Volumetric fluid flux exchanged between branches of a intersection \\
    $\mu$ & Parameters \\
    $\mu_{geom}$ & Geometrical parameters \\
    $\mu_{phy}$ & Physical parameters \\
    $\nu$ & Normal of a sliding surface \\ 
    $\rho$ & Activation function \\
    $\sigma$ & Right-hand side of mesh deformation system \\
    $\Tilde{\sigma}$ & Standard deviation \\
    $\theta$ & Scalar function $\mu$-dependent \\
    $\Theta$ & Parameter space \\
    $\Sigma$ & Singular values matrix \\
    $\upsilon$ & Unit vector \\
    $\upsilon_\gamma$ & Unit vector associated to $\gamma$ \\
    $\hat{\upsilon}$ & Unit vector aligned with the fault \\
    $\varphi$ & Map $\varphi : \Theta \rightarrow \mathcal{V}_n$. In the DL-ROM approach, it is represented by the reduced map network \\
    $\Phi$ & Transition matrix \\
    $\Psi$ & Map $\Psi : \mathcal{V}_n \rightarrow \mathcal{V}_N$. In the DL-ROM approach, it is represented by a decoder \\
    $\Psi'$ & Map $\Psi' : \mathcal{S} \rightarrow \mathcal{V}_n$. In the DL-ROM approach, it is represented by an encoder \\
    $\Omega$ & Matrix domain \\
    $\partial \Omega$ & Boundary of $\Omega$ \\
    $\partial_p \Omega$ & Boundary of $\Omega$ where Dirichlet boundary condition for the pressure is applied \\
    $\partial_q \Omega$ & Boundary of $\Omega$ where Neumann boundary condition is applied \\
    $\partial_{ex} \Omega$ & External boundary of $\Omega$ \\
    $\partial_{in} \Omega$ & Internal boundary of $\Omega$ \\
\end{longtable}

\end{document}